\magnification  1200
\ifx\eplain\undefined \input eplain \fi

\baselineskip13pt 
\vsize=9.1truein

\overfullrule=0pt


\def\wtilde{\widetilde}

\font\smalsmalbf=cmbx8
\font\tenbi=cmmib10

\font\eightbi=cmmib9
\font\sixbi=cmmib6
\font\fivebi=cmmib5
\newfam\bmifam\textfont\bmifam=\tenbi\scriptfont
\bmifam=\eightbi\scriptscriptfont\bmifam=\fivebi

\font\smalltenrm=cmr8
\font\smallteni=cmmi8
\font\smalltensy=cmsy8
\font\smallsevrm=cmr6   \font\smallfivrm=cmr5
\font\smallsevi=cmmi6   \font\smallfivi=cmmi5
\font\smallsevsy=cmsy6  \font\smallfivsy=cmsy5
\font\smallsl=cmsl8      \font\smallit=cmti8

\def\smallfonts{\lineadj{80}\textfont0=\smalltenrm  \scriptfont0=\smallsevrm
                \scriptscriptfont0=\smallfivrm
    \textfont1=\smallteni  \scriptfont1=\smallsevi
                \scriptscriptfont0=\smallfivi
     \textfont2=\smalltensy  \scriptfont2=\smallsevsy
                \scriptscriptfont2=\smallfivsy
      \let\it\smallit\let\sl\smallsl\smalltenrm}

\font\eightbi=cmmib10 at 8pt
\font\sixbi=cmmib10 at 7pt
\font\smathbold=msbm8
\font\mathbold=msbm9 at 10pt\def\Sn{{S^{n-1}}}

\def\N{{\hbox{\mathbold\char78}}}\def\S{{\cal S} }
\def\R{{\hbox{\mathbold\char82}}}
\def\sR{{\hbox{\smathbold\char82}}}

\def\ssRn{{{\hbox{\sixbi\char82}}^{\smallfivrm\char110}}}

\def\avg{{-\hskip-1em\int}}

\def\e{\epsilon}

\def\P{{\cal P}}

\def\d{\delta}


\def\lineadj#1{\normalbaselines\multiply\lineskip#1\divide\lineskip100
\multiply\baselineskip#1\divide\baselineskip100
\multiply\lineskiplimit#1\divide\lineskiplimit100}

\def\remark#1.{\medskip{\noin\bf Remark #1.\enspace}}
\def\endpf{$$\eqno/\!/\!/$$}

\def\pf#1.{\smallskip\noin{\bf  #1.\enspace}}

\def\noin{\noindent}

\def\ds{\displaystyle}
\def\ts{\textstyle}

\def\e{\epsilon}\def\part{\partial_t}

\def\isn{\int_{S^{n-1}}}\def\wtilde{\widetilde}
\def\what{\widehat}
\def\Rn{\R^n}

\def\irn{\int_{\sR^n}}
\def\ref#1{{\bf{[#1]}}}

\def\p{\partial}

\def\half{{\ts{1\over2}}}

\def \L{{\cal L}}

\def\U{{\cal U}}

\def\s{\sigma}
\def\supp{{\rm supp }\,}
\def\b{\beta}
\def\a{\alpha}
\def\na{{n/\a}}\def\nsa{{n\over \a}}\def\an{{\a/n}}\def\asn{{\a\over n}}
\def\nna{{n\over n-\a}}
\def\bp{{\b'}}\bigskip\def\b{\beta}\def\bi{{1\over\b}}\def\bpi{{1\over\b'}}
\def\nf{{\|f\|_\na^\na}}\def\ntf{{\|T f\|_\na^\na}}
\centerline{\bf  Sharp exponential integrability for  critical Riesz potentials}
\centerline{\bf  and fractional Laplacians   on $\R^n$}
\bigskip \centerline{Luigi Fontana, Carlo Morpurgo}
\footnote{}{\smallfonts
\hskip-2.4em  This work was partially supported by NSF Grant DMS-1401035 and by Simons Foundation Collaboration Grant 279735}
\midinsert
{\smalsmalbf Abstract. }{\smallfonts We derive  sharp Adams inequalities for the Riesz and more general Riesz-like  potentials on the whole of $\sR^n$.
As a consequence, we obtain sharp Moser-Trudinger inequalities for the critical Sobolev spaces $W^{\a,\nsa}(\sR^n)$, $0<\a<n$. These inequalities involve fractional Laplacians, higher order gradients, general homogeneous elliptic operators with constant coefficients, and  general trace type Borel measures.
}

\endinsert\medskip

\centerline{\bf 1. Introduction and main results}\bigskip
A classical result in Analysis states that the Riesz potential 
$$I_\a f (x)=\irn|x-y|^{\a-n}f(y)dy,\qquad 0<\a<n\eqdef{rieszpotential}$$
maps $L^p(\Rn)$ continuously onto $L^{np\over n-\a p}(\Rn)$, for $1<p<\nsa$, in other words
we have the following Sobolev inequality 
$$\|I_\a f\|_q\le A\|f\|_p,\qquad q={np\over n-\a p},\qquad f\in L^p(\R^n).\eqdef{HLS}$$ 
For $p=2$ the sharp constant $A$ has been found by Lieb, in his celebrated paper [L].

The critical case $p=\nsa$ has been thoroughly analyzed on domains with finite measure. Indeed, in [A] Adams showed that  if $E\subseteq \R^n$  has finite measure then there is $C>0$, depending on $\a,n,|E|$,  such that for every $f\in L^{\nsa}(\R^n)$ with $\supp f\subseteq E$ and $\|f\|_\na\le1$
$$\int_E \exp\bigg[{1\over|B_1|}|I_\a f(x)|^{\nna}\bigg] dx\le C,\qquad |B_1|={\omega_{n-1}\over n}\eqdef{*2}$$
where $\omega_{n-1}$ is the volume of the $(n-1)-$dimensional sphere, and where uniformity in $f$ is lost if the exponential constant $|B_1|^{-1}$ is replaced by a larger constant. 

In light of this result, it is natural to wonder whether an inequality such as \eqref{*2} can hold without any restrictions on the support of $f$.
A small refinement of Adams' proof (see [FM4], Thm. 3 and Thm. 6) yields the inequality
$$\int_E \exp\bigg[{1\over|B_1|}|I_\a f(x)|^{\nna}\bigg] dx\le C(|E|+|\supp f|)\eqdef{2**}$$
for all compactly supported $f$ with $\|f\|_\na\le 1$, where $C$ depends only on $\a$ and $n$.
However, for any given $E$ with positive measure, one cannot hope to have the exponential integral be uniformly bounded under the sole condition that $\|f\|_\na\le 1$, even if the exponential constant is replaced by any arbitrary positive number. One can easily see this 
along the family of  functions $\chi_{2\le |x|\le R}^{}|x|^{-\a}/\log|x|$, or by a suitable dilation argument. 

Our first main result of this paper is that the sharp Adams inequality holds uniformly with respect to $|\supp f|$ if further conditions on $\|I_\a f\|_p$ are imposed:

\proclaim Theorem 1. For $0<\a<n$, there exists $C=C(\a,n)$ such that for all compactly supported $f\in L^\nsa(\Rn)$ with
$$\|f\|_\na^\na+\|I_\a f\|_\na^\na\le 1\eqdef{1?}$$
and for all measurable $E\subseteq \R^n$ with $|E|<\infty$, we have
$$\int_E \exp\bigg[{1\over|B_1|}|I_\a f(x)|^{\nna}\bigg] dx\le C(1+|E|).\eqdef{2?}$$
Moreover, if $|E|>0$ the exponential constant in \eqref{2?} is sharp, that is, it cannot be replaced by a larger constant.\par

 Note that the potential  of a function in  $L^\nsa(\Rn)$ as defined in \eqref{rieszpotential} need not be even pointwise finite a.e.; for example, when $f(x)=\chi_{|x|\ge2}^{}|x|^{-\a}/\log|x|$ the integral in \eqref{rieszpotential} is everywhere infinite. Using truncates of the same functions it is also clear that  $I_\a$ cannot be extended to a continuous operator on $L^\nsa(\Rn)$. 
This fact is indeed at the origin of the   main difficulties encountered in the proof of Theorem 1. In general the values of  $I_\a f$ strongly depends on the uncontrolled support of $f$. However, by using \eqref{1?} we are able to neutralize this lack of control through a series of careful decompositions of $E$, $f$, and $I_\a f$.

We should point out that it is quite possible that for some $\Omega\subseteq \Rn$ with infinite measure and for some $\alpha\in(0,n)$,  the mapping $f\rightarrow  (I_\a f)|_{\Omega} $ defines  a continuous operator  from $L^\nsa(\Omega)$ to itself,  where we identify $L^p(\Omega)$ with $\{f\in L^p(\Rn): \;f=0$ a.e. on $ \Omega^c\}$. In such cases estimate \eqref{2?} holds under the weaker condition $\|f\|_\na\le 1$ (this point will be addressed in our forthcoming paper [FM3].)

\bigskip

The second natural question to ask is whether Theorem 1 can be extended to a larger space of functions in $L^\nsa(\Rn)$. 
 The operator $I_\a$ is not a closed operator as defined in the space
$$D_0(I_\a):=\{f\in L_c^\na(\Rn):\; I_\a f \in L^\nsa(\Rn)\}$$
where we denoted
$$L_c^p(\Rn)=\{f\in L^p(\Rn):\,\,\supp \!f \;{\rm { compact  }}\}.$$

If $f\in L_c^\na(\Rn)$ then $I_\a f\in L^\nsa(B)$, for any given ball $B$ (this follows for example from  O'Neil's lemma, see Lemma 13 below). 
However, we cannot expect to have $I_\a f$ to  be in $L^p$ outside a large ball, unless $f$ has enough vanishing moments.  
To clarify this, note that   $I_\a f(x)\sim |x|^{\a-n}\int_\ssRn f$, as $|x|\to\infty$. This implies that $I_\a f\in L^\nsa(\Rn)$ if  $0<\a<{n\over2}$, in which case $D_0(I_\a)=L_c^\na(\Rn)$.  
One easily sees, via the Taylor expansion of $|x|^{\a-n}$, that  if ${n\over2}\le \a<{n+m\over2}$, $m=1,2,...n$, then $D_0(I_\a)$ contains the class of $f\in L_c^\na(\Rn)$ with vanishing moments up to order $m-1$. For low values of $m$ one can check that  these spaces actually coincide, and they likely coincide for all values of $m$ (see Remarks after Lemma 15 in Section 6.)

In Theorem 7, we will prove that $I_\a$ is closable, and its smallest closed extension, still denoted as $I_\a$,  will have a domain $D(I_\a)$,  the closure of $D_0(I_\a)$ under the norm in \eqref{1?}, or any other equivalent norm.
 As a consequence, we are able to extend the validity of Theorem 1 to all functions in this larger domain $D(I_\a)$. In the same Theorem we will prove that the image of $D(I_\a)$ under $I_\a$ is indeed the entire Bessel potential space $W^{\a,\nsa}(\Rn).$

 Theorem 1 can be formulated  for sets $E$ having infinite measure, in particular $E=\Rn.$ In this case the exponential needs to be regularized, otherwise the integral in \eqref{2?} is trivially infinite. The natural way to do this is to consider 
$$\exp_N(t)=e^t-\sum_{k=0}^N{t^k\over k!},\qquad N=0,1,...$$
in particular for the value $N=[\nsa-2]$, where $[x]$ denotes the ceiling of $x$, i.e. the smallest integer greater or equal $x$, for $x\in \R$. This type of  regularization has been considered by several  authors in the context of Moser-Trudinger inequalities on spaces of infinite measure (see discussion and references below). Note that if $ {n\over2}\le\a<n$, then $\exp_{[\nsa-2]}^{}(t)=e^t-1$.

\proclaim Corollary 2. If $\big(I_\a,D(I_\a)\big)$ denotes the smallest closed extension of $\big(I_\a,D_0(I_\a)\big)$, then there exists $C=C(\a,n)$ such that for all $f\in D(I_\a)$ satisfying $\|f\|_\na^\na+\|I_\a f\|_\na^\na\le 1$ inequality \eqref{2?} holds,  for all $E$ with finite measure. Moreover, for all such $f$ we also have  
$$\irn \exp_{[\nsa-2]}\Big({1\over|B_1|}|I_\a f(x)|^{\nna}\Big)dx\le C,\eqdef{2??}$$
and the exponential constant in \eqref{2??} is sharp.\par

It is worthwhile noting that the inequality for the regularized exponential \eqref{2??} is an elementary consequence of the one over sets of finite measure given in \eqref{2?}. The proof of this fact is rather straightforward: write $\R^n=E\cup E^c$ with $E=\{x: |I_\a f (x)|\ge1 \}$, and split the integral in \eqref {2??} accordingly. Given  $\|I_\a f\|_\na\le1$ we have  $|E|\le 1$,  and the integral over $E^c$ can be estimated by $\|I_\a f\|_\na$, writing the exponential as a  Taylor series. The real challenges then  are in the region where the potential is large, e.g. in the set 
$\{|I_\a f(x)|\ge1\}$. This observation applies also to  all the other similar inequalities involving the exponential integral on $\Rn$ or on spaces with infinite measure. In short, the regularized exponential is nothing more than a gimmick, and can be replaced by the usual exponential provided that  the resulting  inequality holds on measurable sets of measure no greater than~1. In section 2  we will state an elementary ``Exponential Regularization Lemma", which will be used implicitly throughout this paper, in order to pass from inequalities over sets of finite measure to inequalities over the whole space. 

Theorem 1 was our main motivation, however in this paper we will obtain  more general results at a minimal extra cost. Specifically, we will obtain a version of Theorem~1 and Corollary 2 for convolution operators of the form 
$$T f(x)=\int_{\sR^n}K(x-y) f(y) dy\eqdef{104}$$
where $K(x)$ is what we will call {\it Riesz-like kernel}, namely a mildly regular  function \break$K:\Rn\setminus\! 0\to\R$  that behaves as a Riesz kernel near 0 and that it is globally estimated by the Riesz kernel:
$$ K(x)=g(x^*)|x|^{\a-n}+O(|x|^{\a-n+\d}),\qquad x^*={x\over|x|},\quad\d>0,\, \quad g:\Sn\to\R\eqdef{riesz1}$$
$$|K(x)|\le C|x|^{\a-n}\qquad x\neq 0.\eqdef{riesz2}$$
 The sharp exponential constant in this case is given by 
$${n\over{\isn|g(\omega)|^{\nna}d\omega}}\eqdef{105}$$
provided that $K$ is smooth enough. This result will be also derived for vector-valued  kernels $K(x)$, whose components are Riesz-like; in this case the operator in \eqref{104} will be acting on vector-valued functions, and the sharp constant is still given as in \eqref{105} (see   Theorem 5).

In the second part of the paper we will use the above results to obtain a number of new sharp Moser-Trudinger inequalities on the whole $\Rn$, for higher order gradients, arbitrary powers of the Laplacian, and even general homogeneous elliptic differential operators with constant coefficients. To describe these results let us recall some basic notation, terminology, and facts.

\def\am{{\a\over2}}
For $\a$ a positive  integer and $p\ge1$, the Sobolev space $W^{\a,p}(\Rn)$ is the closure of $C_c^\infty(\Rn)$ under the norm $\sum_{|k|\le \a}\|D^k u\|_p$, where $k$ is a multi-index and $D^ku$ denotes the $k$-th distributional derivative of $u$. Alternatively, $W^{\a,p}(\Rn)$ is the space of $u\in L^p$ such that $D^k u \in L^p(\Rn)$.  It is well-known that when $p>1$ the Sobolev space can be characterized in terms of Bessel potentials as
$$W^{\a,p}(\Rn)=\{u\in \S':\;(I-\Delta)^\am u\in L^p(\Rn)\}=\{G_\a*f,\,f\in L^p(\Rn)\}\eqdef{W}$$
with equivalent norm $\|(I-\Delta)^\am u\|_p$, where $\S$ is the space of Schwarz functions, $\S'$ is the space of tempered distributions,  and  $G_\a$ is the Bessel potential, i.e. the $L^1$ function such that  $\what G_\a(\xi)=(1+4\pi^2|\xi|^2)^{-\am}$.  When $\a$ is an arbitrary positive real number
and $p>1$ we will still use the notation in \eqref{W} and call $W^{\a,p}(\Rn)$  Sobolev space, with the understanding that it is really  the Bessel potential space, and it is not to be confused with the Aronszajn-Gagliardo-Slobodeckij space.  In the literature various other notations exist for such spaces ($H^{\a,p},\, \L_\a^p$, etc.).

\def\Da{(-\Delta)^{\am}}

If $\Delta$ denotes the Laplace operator, we define $\Da \phi$ for arbitrary $\a>0 $ as the inverse Fourier transform  of $(2\pi|x|^\a)\what\phi$, for any $\phi\in \S$. This defines a $C^\infty$ function satisfying the estimate
$$|\Da \phi(x)|\le C{p(\phi)\over(1+|x|)^{\a+n}},
\eqdef{schwarz}$$
where $p(\phi)$ is a suitable finite sum of the usual seminorms of $\phi$ in $\S$, and $C$ is independent of $\phi$.  Such an estimate is proven for example in [GO, Lemma 1], even though the dependence of the constant on the seminorms in $\S$, is not explicitly written there  (but it can be easily  derived from the proof.) See also [Hy2, Lemma 2.2]. With the aid of \eqref{schwarz} one can  extend the definition of $\Da u$ as a tempered distribution for any $u\in L^p$, $1\le p<\infty$.

 It is an exercise in distribution theory (see Section 7 for some more details) to show that for $\a>0$ and $p>1$
$$W^{\a,p}(\Rn)=\{u\in L^p(\R^n):\;\Da u\in L^p(\Rn)\}$$
with equivalent norm $\|u\|_p+\|\Da u\|_p$, or any norm of type $\big(\|u\|_p^{q}+\|\Da u\|_p^{q}\big)^{1\over q}$, $\;1\le q\le\infty$. In Theorem 7 we will prove that when $0<\a<n$ the space of critical Riesz potentials indeed coincides with the whole $W^{\a,\nsa}(\Rn)$ i.e.
$$W^{\a,\nsa}(\Rn)=\{I_\a f,\,f\in D(I_\a)\},\eqdef{sob}$$
and in this space the operator $\Da$ is invertible, with  inverse given by the integral operators having kernel $c_\a I_\a$, where
$$ c_\a={\Gamma\big({n- \a\over2}\big)\over 2^\a\pi^{n/2}\Gamma\big({ \a\over2}\big)}.\eqdef{calpha}$$
This characterization of the critical Sobolev space, which can be likely extended to any $W^{\a,p}$, with $p\ge \nsa$,  does not seem to have been considered before. Samko [Sa, Thm. 7.18]  proves an even  more general result, but with a  $I_\a f$ defined as a Lizorkin distribution on  $L^p(\Rn)$, for any $p\ge \nsa$. It can be seen that for $f\in D(I_\a)$ our definition of $I_\a f$ coincides with the one given by Samko, however our approach appears to be more  direct and useful   in the context of exponential integrability.

Finally, recall that  a homogeneous  elliptic differential operator of  even order $\alpha<n$ with real constant coefficients has form
$$Pu=\sum_{|k|=\a}a_k D^k u\eqdef{P}$$
 acting, say, on $C^\infty_c(\Rn)$, with 
$$p_\a(\xi)=P(2\pi i \xi)=(2\pi)^\a(-1)^{\a/2}\sum_{|k|=\a}a_k \xi^k,\quad |p_\a(\xi)|\ge c_0|\xi|^\a,\qquad \xi\in \R^n$$
for some $c_0>0$. The fundamental solution of $P$ is given by a convolution operator with kernel $g_P^{}$ given by 
$${g_P^{}}^{}(x)=\irn {e^{-2\pi i x\cdot \xi}\over p_\a(\xi)}d\xi\eqdef{gp}$$
in the sense of distributions.
\def\Dao{\nabla(-\Delta)^{{\a-1\over2}}}
\proclaim Theorem 3. For $0<\a<n$ let  $P$ be either $\Da$, $\Dao$ for $\a$ odd, or a homogeneous elliptic operator of even order $\a<n$ with constant coefficients. Then there exists $C=C(\a,n,P)>0$ such that for every $u\in W^{\a,\nsa}(\Rn)$ with 
$$\|u\|_\na^\na+\|Pu\|_\na^\na\le 1\eqdef{K5}$$
 and for all measurable $E\subseteq \Rn$ with $|E|<\infty$ we have
$$\int_E\exp\bigg[\gamma(P)|u(x)|^\nna\bigg]dx\le C(1+|E|)\eqdef{K4}$$
and
$$\int_{\sR^n}\exp_{[\nsa-2]}\bigg[\gamma(P)|u(x)|^\nna\bigg]dx\le C\eqdef{K4a}$$
where
$$\gamma(P)=\cases{\ds{c_\a^{-\nna}\over|B_1|}, & if $P=\Da$, any $\a\in(0,n)$    \cr\ds{\big((n-\a-1)c_{\a+1}\big)^{-{\nna}}\over|B_1|} &  if $P=\Dao$ and $\alpha$ odd.\cr}\eqdef{gamma1}$$
and where
$$\gamma(P)=\ds{n\over\isn |{g_P^{}}^{}(\omega)|^\nna d\omega},\eqdef{gamma2}$$
 if $P$ elliptic as in \eqref{P} and $\a$ even.
Moreover, the exponential constant $\gamma(P)$ in \eqref{K4} and  \eqref{K4a} is sharp. \par
For a more general sharpness statement see \eqref{sharper} and corresponding remarks.

Theorem 3 in the present form is only known so far in the case $\a=1$, and $\alpha=2$. The earliest result is for $P=\nabla$ and it  is due due to Ruf [Ruf], who proved it in dimension~2. Ruf's result was  later extended to all dimensions by Li-Ruf [LR]. The case $P=\Delta$ was settled by Lam-Lu in [LL2, Thm 1.5], who reduced the problem to the sharp Moser-Trudinger inequality on bounded domains with homogeneous Navier boundary conditions derived by Tarsi [Tar]. The case   $P=(-\Delta)^{\half}$ in dimension $n=1$ was settled recently  by Iula-Maalaoui-Martinazzi [IMM]. We note here that two questions posed in [IMM, Sec. 1.3] are positively answered by our Theorem 3 (see remarks after the proof of Theorem 3, Section 8). 

For $P=\nabla$ and $P=\Delta$ there are in fact even more refined inequalities due to Ibrahim-Masmoudi-Nakanishi [IMN] ($n=2$, $P=\nabla$),  Masmoudi-Sani [MS1] ($n\ge 2$, $P=\nabla$), [MS2] ($n=4$, $P=\Delta$), and Lu-Tang-Zhu [LTZ] ($n\ge3$, $P=\Delta$),  who proved that  
$$\int_{\sR^n}{\exp_{[\nsa-2]}\Big[\gamma(P)|u(x)|^\nna\Big]\over(1+|u(x)|)^{\nna}}dx\le C, \eqdef{masmoudi}$$
under the weaker norm condition 
$$\max\big\{\|u\|_{\na},\,\|P u\|_\na\big\}\le 1,\qquad u\in W^{\a,\nsa}(\Rn).\eqdef{Adachi1}$$
It is possible to see that this result is stronger than our Theorem 3, however it is only known  for $P=\nabla$ and $P=\Delta$. 
\footnote{(*)}{\smallfonts
  See ``Note added in proof'' at the end of the paper}

Rearrangement tools such as the P\'olya-Szeg\"o inequality or Talenti's inequality  were crucial ingredients in each  of the  known results cited above,  but are completely avoided in our paper, which is based on exponential integrability of potentials, or Adams type inequalities.

 In the context of the Heisenberg group the Li-Ruf result was obtained by Lam and Lu [LL1], and on general noncompact Riemannian manifolds by  Yunyan Yang [Y], but only in the subcritical case, i.e. for exponential constants $\gamma<\gamma(\nabla)$.

Inequality \eqref{K4a} under the weaker norm condition \eqref{Adachi1}
was first obtained without sharp exponential constants by  Ogawa [Og], for $P=\nabla$ and  $n=2$. Ogawa's result was later extended to $P=\Da$, for any $n$ and any $\a\in (0,n)$ by Ozawa. 

When $P=\nabla$ inequality \eqref{K4a} was derived for any exponential constant $\gamma<\gamma(\nabla)$ first by Cao [Cao] in dimension 2, followed by Panda [Pa] and Do \'O [Do\'O] in any dimension. A couple of years later Adachi-Tanaka [AT] reproved the same result and cast it in a dilation invariant form
$$\int_{\sR^n} \exp_{n-2}\bigg[\gamma\bigg({|u(x)|\over \|\nabla u\|_n}\bigg)^{n\over n-1}\bigg]dx\le C\bigg({\|u\|_n\over\|\nabla u\|_n}\bigg)^n,\qquad u\in W^{1,n}(\Rn)\setminus\{0\}\eqdef{*7}$$
where $\gamma<\gamma(\nabla)$.
 Moreover, Adachi and Tanaka [AT] verified  that the above inequality \eqref{*7} fails if $\gamma=\gamma(\nabla)$ by means of the usual Moser sequence, along which the left-hand side is bounded away from 0 (in fact bounded above too by Moser's inequality), whereas the  right-hand side tends to~0. 

As a corollary of Theorem 3 we have an Adachi-Tanaka type result, with sharp control on the right-hand side constant, for  any of the equivalent  norms
$$\|u\|_{\na,q}:=\cases{\big(\|u\|_\na^{q\na}+\|Pu\|_\na^{q\na}\big)^{\a\over qn} & if  $\;1< q<\infty$\cr\cr
\max\big\{\|u\|_\na,\|Pu\|_\na\big\} & if $\;q=+\infty,$\cr}\eqdef{qnorm}$$
each of which yields a norm in the Sobolev space $W^{\a,\nsa}(\Rn)$, even for $q\ge{\a\over n}$. Observe that $\|u\|_{\na,q}^q$ is decreasing in $q$ when $\|u\|_{\na,q}\le 1$, so that the conclusions of Theorem~3 are true  if condition \eqref{K5} is replaced by $\|u\|_{\na,q}\le 1$, for all $q\in\big[{\a\over n},1]$. The following  result shows that this is no longer the case for $q>1$:
\eject
\proclaim Corollary 4. For $0<\a<n$ and  $P$  as in Theorem 3, there exists $C=C(\a,n,P)$ such that for any $\theta\in (0,1),\; 1<q\le+\infty$,  for every  $u\in W^{\a,\nsa}(\Rn)$ with $$\|u\|_{\na,q}\le 1,\eqdef{103ab}$$ and   for all  measurable  $E\subseteq \Rn$ with $|E|<\infty$   we have
$$\int_E \exp\Big[\theta\gamma(P)|u(x)|^\nna\Big]dx\le C( 1-\theta)^{-{1\over q'}}(1+|E|),\eqdef{102ab}$$
and
$$\int_{\sR^n}\exp_{[\nsa-2]}\Big[\theta\gamma(P)|u(x)|^\nna\Big]dx\le C( 1-\theta)^{-{1\over q'}}.\eqdef{102bc}$$
The inequalities \eqref{102ab},\eqref{102bc} are sharp in $\theta$, in the sense that for given $q\in(1,\infty]$,  there exists a family of functions $\{u_\theta\}\in W^{\a,\nsa}(\Rn)$ satisfying \eqref{103ab} for which  \eqref{102ab} (for any given $E$ with positive measure, with $C$ depending on $E$) and \eqref{102bc} are reversed; in particular, the exponential integral cannot be uniformly bounded if $\theta=1$.  
\par

See also  the recent  papers [CST] and  [LLZ], where the relations between the Li-Ruf and the Adachi-Tanaka type results as in Corollary 4  are explored in more detail. 

 Corollary 4 suggests that among all possible natural Sobolev norms  in $W^{\a,\nsa}(\Rn)$ which involve  only $u$ and $P u$,  the ``Ruf norm" $\|u\|_{\na,1}$  yields  the least restrictive condition under which the  Moser-Trudinger inequality of the type \eqref{K4a} holds. Under norm conditions more restrictive than \eqref{K5} sharp higher order results do exist. Indeed, if $P=(I-\Delta)^{\a\over2}$, then the  sharp inequality in \eqref{K4a} holds  for any $\alpha\in(0,n)$ under the condition  
$$\|(I-\Delta)^{\a\over2}u\|_\na\le 1.\eqdef{*15}$$
This result goes back to Adams, who proved it for $\alpha=2$  in his original 1988 paper [A, Thm. 3]. Strictly speaking, Adams proved that if  $\Omega$ is an open and bounded set with measure $|\Omega|\le 1$,  $\alpha=2$ and if $u\in W^{2,{n\over2}}(\R^n)$ satisfies \eqref{*15}, then the basic Moser-Trudinger inequality \eqref{K4a} holds with exponential constant  $\gamma(\Delta)$. To prove this result Adams  modified slightly the proof that he gave of \eqref{*2} for Riesz potentials on   $L^\nsa(\Omega)$, adapting it   to Bessel potentials on $L^{\nsa}(\R^n)$; in particular, he proved that for $\alpha=2$ there exists $C>0$, independent of $\Omega$, such that 
$$\int_\Omega \exp\bigg[{1\over|B_1]} |G_\a* f(x)|^{\nna}\bigg] dx\le C,\qquad \|f\|_{\na}\le 1\eqdef{*16}$$
where $G_\alpha$ denotes the usual Bessel potential. In this setting it is not necessary to have $\Omega$ be a bounded open set, it's enough that $\Omega$ be measurable with measure $\le 1$. Adams' proof of \eqref{*16} for $\alpha=2$, however,  is actually working for any $\alpha\in (0,n)$, after straightforward modifications; it appears that the only reason why Adams considered $\alpha=2$ is because in this case the norm in \eqref{*15} coincides precisely with the usual full Sobolev norm on $W^{2,{n\over2}}$. The fact that the Bessel potential kernel $G_\a$ decays quite well  at infinity is what, 
at the end of the day, makes Adams' proof go through with few modifications from the finite measure case. These issues are thoroughly addressed in our paper  [FM3].

In the same spirit as in Corollary 4, inequality  \eqref{*16}, valid for any measurable $\Omega$ with $|\Omega|\le1$,  implies easily that (and it is in fact equivalent to) 
$$\int_{\sR^n} \exp_{[\nsa-2]}\Big({1\over|B_1|}|G_\a*f|^{\nna}\Big)dx\le C,\qquad \|f\|_{\na}\le 1,\eqdef{*17}$$
however Adams does not mention this form of the inequality in his paper (for the case $\a=2$).

Unaware of Adams' result, Ruf-Sani [RS] proved \eqref{K4a} under the condition \eqref{*15} when $\alpha$ is an even integer,  not using the Bessel potential approach, but rather comparison theorems for  suitable Navier boundary value problems.
Recently, Lam-Lu [LL2]  presented a complete proof of the same result for any $\alpha\in (0,n)$, using the Bessel potential and following Adams' original argument.  

In [FM3] we will present a general Adams inequality on measure spaces, extending our original results in [FM1] to integral operators over sets of arbitrary measure, whose kernel satisfies additional decay conditions at infinity.   As a special case, an inequality like \eqref{*17}  holds for convolution operators $K*f$ on $\Rn$, under the conditions that $K (x)\sim |x|^{\a-n}$ around 0, and $K\in L^{\nna}\cap L^\infty$ at infinity. This result complements nicely the Adachi-Tanaka type result for Riesz-like potentials, in Theorem 6 below, where we show that the Adams inequality fails at the critical exponential constant if $K$ is Riesz-like and {\it not} in $L^\nna$ at infinity.

Another perhaps interesting observation, is that Corollary 4 implies the following weak Masmoudi-Sani type inequalites:

$$\int_{\sR^n}{\exp_{[\nsa-2]}\Big[\gamma(P)|u(x)|^\nna\Big]\over(1+|u(x)|)^{\nna}}dx\le Cq, \eqdef{masmoudi}$$
valid for $1\le q<\infty$ and for all $u$ such that $\|u\|_\na^{qn/\a}+ \|\nabla^\a u\|_\na^{q\na}\le 1$, and
$$\int_{\sR^n}{\exp_{[\nsa-2]}\Big[\gamma(P)|u(x)|^\nna\Big]\over(1+|u(x)|)^{\nna(1+\e)}}dx\le {C\over\e}, \eqdef{masmoudi2}$$
for all $\e>0$, and for  $u$ such that $\|u\|_\na\le 1$ and $\|\nabla^\a u\|_\na\le 1$. These inequalities can be obtained by taking $E=\{|u|\ge 1\}$ in \eqref{102}, and integrating in $\theta$ directly, or after multiplying by  $(1-\theta)^\e$. Better inequalities can be obtained with this method, but still without reaching the optimal Masmoudi-Sani type of result.

Our method also allows us to obtain, with minimal modifications,  the results in Theorems 1 and  3 when the non-regularized exponential is integrated against a general Borel measure $\nu$ satisfying
$$\nu\big(B(x,r)\big)\le Q r^{\sigma n},\qquad \forall x\in \Rn,\;\;r>0\eqdef{meas}$$
for some $Q>0$ and $\sigma\in (0,1]$.
The resulting Adams/Moser-Trudinger inequalities analogous to \eqref{2?} and \eqref{K4} are of trace type, and the sharp constant in this case is $\sigma \gamma(P)$. Similar results hold for the general operators as in \eqref{104}. The first results regarding Moser-Trudinger inequalities of trace type for general measures as in \eqref{meas} on bounded open sets, are due to Cianchi [Ci1], for the case $\alpha=1$. In [FM1] the authors extended Cianchi's results to higher order gradients.  Examples of measures as in \eqref{meas} are the Hausdorff measures on  submanifolds of  $\Rn$. Another example  is the ``singular measure" with density $|x|^{(\sigma-1)n}$, which was considered for example in [LL2] and [AY], and in other papers dealing with domains of finite measure.  For those measures, and even for more general ones,  we can obtain the results in the full regularized form over the whole $\Rn$, as in \eqref{2??}, \eqref{K4a}. For sake of clarity we have not stated our main results for such more general measure, however  in  Theorem 18  we will provide an  explicit statement.

A byproduct of the technique used in the proof of Theorem 1 is the following sharp Trudinger inequality without boundary conditions on bounded smooth domains $\Omega$:
$$\int_\Omega \exp\bigg[2^{-{1\over n-1}}\gamma(\nabla)|u(x)|^{n\over n-1}\bigg]dx\le C\eqdef{mtw}$$
 for each $u\in W^{1,n}(\Omega)$ with $\|u\|_n^n+\|\nabla u\|_n^n\le 1$ (see Theorem 19). This result appears to be the only known sharp version of the original inequality due to Trudinger [Tr, Thm. 2].  Under the condition $\|\nabla u\|_n\le 1$ and $u$ with 0 mean, the sharp constant was found in [CY] when $n=2$ and in [Ci2] for general $n$; a version of their result for $u\in W^{2,n/2}$ in the unit ball was given by the authors in [FM2]. 

Finally, we believe that the techniques developed in this paper could be adapted  to other settings, such as the Heisenberg group or other noncompact manifolds. For example, on the Heisenberg group one could consider a version of Theorem 3 for the powers of the sublaplacian, which would extend the results in [LL1] to higher order operators.

\vskip1em
For the convenience of the reader  here is  an outline of how the paper is organized. \smallskip
\noin{\underbar{\sl Section 2}}.   We introduce Riesz-like kernels and potentials and state three theorems for them:
\item{-} Theorem 5:  A  sharp Adams inequality under a   Ruf type condition. This theorem contains Theorem 1 as a special case.
\item{-} Theorem 6: A sharp Adams inequality for a  family of Adachi-Tanaka type conditions. 
\item{-} 
Theorem 7:  Closability of Riesz-like potentials and extension of results in Theorems~5,6 to larger function spaces. Characterization of classical Sobolev-Bessel potential spaces, at the critical index, in terms of  Riesz potentials.
\smallskip
\noin{\underbar{\sl Section 3.}}  We prepare the ground for the proof of Theorem 5. Several results are proven including:

\item{-} Proposition 8: a refined version of Adams' inequality for sets of finite measure.
\item{-} Lemma 9: the Exponential Regularization Lemma.
\item{-} Lemma 10: A crucial lemma about exponential integrability under additive perturbations.
\item{-} Lemmas 11-14: A series of lemmas about regularity and  integrability properties of Riesz-like potentials.
\smallskip
\noin{\underbar{\sl Section 4}}.  Proof of the  inequalities in Theorem 5.\smallskip
\noin{\underbar{\sl Section 5}}. Proof of  the inequalities  in Theorem 6.
\smallskip
\noin{\underbar{\sl Section 6}}. Proof of the sharpness statements of Theorems 5 and Theorem 6.
\smallskip
\noin{\underbar{\sl Section 7}}.  Proof of  Theorem 7.\smallskip
\noin{\underbar{\sl Section 8}}. Proof of Theorem 3.\smallskip
\noin{\underbar{\sl Section 9}}. Proof of Corollary 4.\smallskip
\noin{\underbar{\sl Section 10}}. We state and prove an extension of Theorems 5 and 3  for a general class of Borel measures (Theorem 18).  \smallskip
\noin{\underbar{\sl Section 11}}. We state and prove a sharp Trudinger inequality on bounded domains without boundary condition  (Theorem 19).
\vskip 1em

\centerline{\bf 2. Adams inequalities for general homogeneous Riesz-like potentials on $\R^n$}\bigskip

Let us define $K\to \R^n\setminus 0\to \R$ to be a {\it Riesz-like kernel} of order $\a\in(0,n)$ if it satisfies the 
following properties:\smallskip
$$K(x)=g(x^*)|x|^{\a-n}+O(|x|^{\a-n+\d}),\qquad x^*=\ds{x\over|x|},\eqdef{rl1}$$
$$ |K(x)|\le C|x|^{\a-n}\qquad x\neq 0\eqdef{rl2}$$
$$|K(x_1)-K(x_2)|\le C|x_1-x_2|\max\big\{|x_1|^{\a-n-1},|x_2|^{\a-n-1}\big\},\qquad x_1,x_2\neq0\eqdef{rl3}$$
for some $\d>0$, $C>0$, $0<\a<n$, and $g:\Sn \to\R$, not identically 0. The ``big O'' notation in \eqref{rl1} means $|O(|x|^{\a-n+\d})|\le C|x|^{\a-n+\d}$ for all $x\neq 0$ near 0 (hence for all $x\neq0$, due to \eqref{rl2}), and  $C$ will denote a constant that may vary from place to place. It is clear that \eqref{rl3} implies that $K$ is Lipschitz outside any ball $B(0,r)$ and that $g$ is Lipschitz on $\Sn$, since $g(x^*)=\lim_{t\to0} K(tx^*)t^{n-\a}$.   Also, if $K$ is homogenous of order $\a-n$, and Lipschitz on $\Sn$, then it is Riesz-like.

We say that a Riesz-like potential is {\it $m-$regular} (with $m\in \N$) if $K\in C^m(\R^n\setminus 0)$ and
$$|D_h^j K(x)|\le C|x|^{\a-n-j},\qquad x\neq0,\; j=1,...,m$$
where $h$ denotes a multi-index with $|h|= j$, and $D_h^jK$ denotes the $j$-th derivative of $K$ w.r. to $h$.
With some straightforward estimates one checks that \eqref{rl3} is implied by the more natural assumption that $K$ is 1-regular, i.e. differentiable and  $|\nabla K(x)|\le C|x|^{\a-n-1}.$. Obviously the Riesz kernel is an $m-$regular Riesz-like kernel for all $m$.

From now on  $T$ will denote the convolution operator with a Riesz-like kernel of order~ $\a$:
$$T f(x)=K* f(x)=\irn K(x-y)f(y)dy.\eqdef{Tg}$$
We define  $T f$ on vector-valued functions in the same way, with the understanding that 
$$K=(K_1,...,K_m),\quad f=(f_1,...,f_m),\quad Kf=K_1f_1+...+K_mf_m,\quad |f|=(f_1^2+...+f_m^2)^{1/2},\eqdef{11b}$$
where each $K_j$ is a Riesz-like kernel of order $\a$.
Additionally, we will let, for $f=(f_1,..,f_m)$
$$f^+=(f_1^+,....,f_m^+),\qquad f^-=(f_1^-,...,f_m^-)\eqdef{11c}$$
where $f_j^+$ and $f_j^-$ denote the positive and negative parts of $f_j$. We will say that a vector $f$ is nonnegative/nonpositive) if each component of $f$ is nonnegative/nonpositive a.e. The results and their proofs below are valid for the vector-valued case with the above conventions;  we will not distinguish between the scalar case and the vector-valued case, except in a few isolated instances.
\smallskip

The following theorem contains Theorem 1 as a special case:\smallskip
\proclaim   Theorem 5.
If $K$ is a nonnegative or nonpositive Riesz-like kernel of order $\a$ , then there
 exists a constant $C=C(\alpha,g,n)$ such that for every measurable $E\subseteq \R^n$ with $|E|<\infty $ and for all compactly supported $f$ with
$$\|f\|_{\na}^\na+\|T f\|_{n/\a}^{n/\a}\le 1\eqdef {11}$$
we have
$$\int_E\exp\bigg[{1\over A_g}\,|T f(x)|^{n\over n-\a}\bigg]dx\le C(1+|E|).\eqdef {12}$$
where $$A_g={1\over n}\isn |g(\omega)|^\nna d\omega,\eqdef A$$
and also
$$\int_{\sR^n}\exp_{[\nsa-2]}\bigg[{1\over A_g}\,|T f(x)|^{n\over n-\a}\bigg]dx\le C.\eqdef {12aa}$$
If $K$ is $n-$regular and  $|E|>0$, then the exponential constant $A_g^{-1}$ in  \eqref {12}, and also in \eqref{12aa}, is sharp, i.e. it cannot be replaced by a larger number.\hfill\break \indent
 If $K$ changes sign then the above results continue to hold for all compactly supported $f$ satisfying the additional pointwise condition: for all $a\in \R^n$
$$|T f(x)|\le \int_{|y-a|\le 2}|K(x-y)|\,|f(y)|dy+C_1\|T f\|_\na,\qquad 
\;|x-a|\le 1\eqdef{12q}$$
almost everywhere, where $C_1$ is a constant depending only on $\alpha$ and $n$.
\par 
\bigskip\def\s{\sigma}

In the same spirit as Corollary 4 we have the following general Adachi-Tanaka type result:

\proclaim Theorem 6. If $K$ is a Riesz-like kernel, then for any $\theta\in (0,1)$ there exists $C_\theta=C(\theta,\a,K,n)$ such that for  $1<q\le+\infty$,   for all   $E\subseteq  \Rn$ with $|E|<\infty$ 
and for all compactly supported $f$ with
$$\cases{\|f\|_\na^{q\na}+\| T f \|_\na^{q\na}\le 1, &  if $q<\infty$\cr\cr   \max\big\{\| f\|_\na,\,\|T f\|_\na\big\}\le 1, & if  $q=+\infty.$\cr}\eqdef{103}$$ we have
$$\int_E \exp\Big[{\theta\over A_g}|T f(x)|^\nna\Big]dx\le C_\theta(1+|E|)\eqdef{102}$$
and also 
$$\int_{\sR^n}\exp_{[\nsa-2]}\Big[{\theta\over A_g}|T f(x)|^\nna\Big]dx\le C_\theta.\eqdef{102aa}$$
If $K$ is $n-$regular and $K\notin L^\nna(|x|\ge1)$ then inequalities  \eqref{102}  and \eqref{102aa} are sharp, in the sense that if $|E|>0$  the exponential integrals cannot be uniformly bounded if $\theta=1$.  
\smallskip
If $K$ is a nonnegative or nonpositive homogeneous Riesz-like kernel then \eqref{102} and \eqref{102aa} hold with $C_\theta$ replaced by $C( 1-\theta)^{-{1\over q'}}$, some $C=C(\a,n,K)$, and the inequalities are sharp in $\theta$, in the sense that they can be reversed along a family $f_\theta\in L^\nsa_c(\Rn)$ satisfying \eqref{103}.
\par

The proof of Theorem 6, and in particular of the original Adachi-Tanaka estimate ($q=+\infty$), without sharp control on the right-hand side in terms of $\theta$, does not require Theorem~5, and it is much simpler to prove it directly with the methods given in this paper. The homogeneous case instead is a direct consequence of Theorem 5, using a dilation argument.

Both Theorem 5 and Theorem 6 are stated in terms of functions $f$ with compact support and such that both $f$ and $T f$  are in $L^\nsa(\Rn)$, i.e. functions $f$ in the space
$$D_0(T):=\{f\in L_c^{\na}(\R^n):\, Tf \in L^\nsa(\Rn)\}.$$
As in the case of the Riesz potential, the space $D_0(T)$ coincides with $L_c^\na(\Rn)$ for $\a<n/2$, whereas for $\a\ge n/2$ it contains  all functions of $L^{\na}_c(\Rn)$ with vanishing moments up to order $n-1$ (see remarks after Lemma 15 in Section 6). It is always possible to normalize any $f\in L_c^1(B)$ ($B$ any ball)  in such a way that all moments up to any given order $m$ are zero,  by subtracting a suitable polynomial of order $m$ restricted to $B$ (see Section 6.)
In the next theorem we show that $T$ has a smallest closed extension, still denoted $T$, and the validity of Theorem 5 and Theorem 6  extends to all functions in the domain of such an extension. Moreover, we can characterize the space $W^{\a,\nsa}(\Rn)$ in terms of the closed extension of the Riesz potential.

\proclaim Theorem 7. If $K$ is a Riesz-like kernel, then  the operator $T: D_0(T)\to  L^\nsa(\Rn)$ is closable, and  its smallest closed extension (still denoted $T$) has domain
$$D(T)=\big\{f\in L^\nsa(\Rn): \,\exists \{f_k\}\subseteq D_0(T),\,\exists h\in L^\nsa(\Rn)\; {\rm with}  \;f_k\buildrel{ L^{\nsa}}\over \longrightarrow f,\;T f_k\buildrel{ L^{\nsa}}\over \longrightarrow h\big\}\eqdef{D}$$
and $Tf=h$. 
 Theorem 5 and Theorem 6 are still valid if the condition $f\in L_c^\nsa(\Rn)$ is replaced by $f\in D(T)$.
\smallskip
In the case of the Riesz potential we have 
$$W^{\a,\nsa}(\Rn)=\{ I_\a f,\, f\in D(I_\a)\}\eqdef{sob1}$$
and the operator $\Da$ is a bijection between $W^{\a,\nsa}(\Rn)$ and $D(I_\a)$, with inverse $c_\a I_\a $.
\par
Note that  if $1<p<\nsa$ the potential $I_\a f$ is well defined on $L^p$ (and belongs to an $L^q$) and  we  automatically get that if $D(I_\a)=\{f\in L^p: \, I_\a f\in L^p\}$  then $I_\a:D(I_\a)\to L^p$ is closed on $D(I_\a)$,  and similarly for the more general $T$. In this sense, there is no need to consider the closure of $I_\a$ for $p$ subcritical. A version of Theorem 7 is likely true in the case $p>n/\a$, perhaps after suitable modifications of the proof  presented later in this paper. 
\bigskip

\centerline{\bf 3. Proof of   Theorem 5 : Overview and preliminary lemmas}\bigskip

Let us start with the following  slightly more refined version of the usual Adams inequality on sets of finite measure, which include \eqref{2**}
as a special case:

\proclaim Proposition 8.  If $K$ is a Riesz-like kernel,  then there exists a constant $C=C(\alpha,g,K)$  such that  for  given measurable sets $E,F$ with finite measure and for all   $f$ with $\supp f\subseteq F$ and  $\| f\|_\na\le1$ we have 
$$\int_E\exp\bigg[{1\over A_g}\,|T f(x)|^{n\over n-\a}\bigg]dx\le C\big(1+|F|\big)\big(1+\log^+|F|+|E|\big).\eqdef {12h}$$
If in addition $K$ is homogeneous, then \eqref{12h} holds with $C(|E|+|F|)$ on the right-hand side.\par
\pf Proof. The result can be proved in essentially the same original argument  by Adams in the general form given  in [FM1], using O'Neil's inequality and  the Adams-Garsia lemma, but keeping better track of the constants (for more details see the original version of the present paper [FM4], Thm. 3 and Thm. 6). 

If $K$ is homogeneous, the  inequality  follows from the case  $|F|=1$, indeed we have the following dilation properties: 
 if  $f_\lambda(x)=\lambda^\a f(\lambda x)$ then 
$$\|f_\lambda\|_\na=\|f\|_\na,\qquad Tf_\lambda(x/\lambda)=T f(x),\qquad \|T f_\lambda\|_\na^\na=\lambda^{-n}\|T f\|_\na^\na.\eqdef{D1}$$
and 
$$\int_E \exp\Big[{1\over A_g}|T f(x)|^\nna\Big]dx=\lambda^n\int_{E/\lambda} \exp\Big[{1\over A_g}|T (f_\lambda)(x)|^\nna\Big]dx.\eqdef{D2}$$
Hence, if the inequality is known for $|F|=1$ it's enough to take $f_\lambda$ as above with $\lambda=|F|^{1/n}$ to obtain the bound $C(|E|+|F|)$  (note that $\supp f_\lambda\subseteq F/\lambda$ which has measure 1). \endpf

\eject

\noin{\bf Remarks.} 

\noin 1)  For {\it fixed} $E$ and $F$ the sharp exponential constant is not $A_g^{-1}$ in general, and  will depend on the relative geometry of $E$ and $F$. 

\noin 2) For the purposes of this paper all we need of Proposition 8  is that the exponential integral in \eqref{12h} is uniformly bounded by $C(1+|E|)$, if $|F|$ is bounded by a given constant.

As we pointed out earlier in regard to \eqref{2**}, one cannot hope to make the right hand side of  \eqref{12h} to be independent of $|F|$ without  restrictions on $\|T f\|_\na$. 
The  ``Ruf condition" in \eqref{11} is precisely  what is needed in order to compensate for the lack of control on the support of $f$, or rather on its measure.

Let us now state an elementary lemma, which is more like an observation, with the hope of  clarifying the equivalence between exponential inequalities on sets of finite measure and regularized exponential inequalities on sets of arbitrary measure.
\smallskip
\proclaim Lemma 9  (Exponential Regularization Lemma). Let  $(N,\nu)$ be a measure space and $1<p<\infty$, $\a>0$. Then for every $u\in L^p(N)$ we have
$$\int_{\{|u|\ge1\}} e^{\alpha|u|^{p'}}d\nu-e^\alpha\|u\|_p^p\le\int_N\bigg(e^{\a |u|^{p'}}-\sum_{k=0}^{[p-2]}{\a^k |u|^{kp'}\over k!}\bigg)d\nu\le\int_{\{|u|\ge1\}} e^{\alpha|u|^{p'}}d\nu+e^\a\|u\|_p^p.$$
In particular, the functional $\int_N \exp_{[p-2]}^{}\big[\a |u|^{p'}\big]$ is  bounded on a bounded subset $X$ of  $L^p$, if and only if $\int_{\{|u|\ge 1\}}\exp\big[\a |u|^{p'}\big]$ is  bounded  on $X$.
\par
\pf Proof. Recall that $[p-2]$ is the smallest integer greater or equal $p-2$ . To start, write $N=(N\cap\{|u|<1\})\cup (N\cap\{|u|\ge1\})$, and split the middle integral accordingly.  Then we just observe that
$$e^{\a|u|^{p'}}-e^\a|u|^p\le e^{\a |u|^{p'}}-\sum_{k=0}^{[p-2]}{\a^k |u|^{kp'}\over k!}\le e^{\a|u|^{p'}},\qquad {\hbox { if }} |u|\ge 1$$
and
$$0\le e^{\a |u|^{p'}}-\sum_{k=0}^{[p-2]}{\a^k |u|^{kp'}\over k!}=\sum_{k=[p]-1}^\infty{\a^k |u|^{kp'}\over k!}\le  e^\a|u|^p,\qquad {\hbox { if }} |u|\le 1.$$

\endpf
Note also that for any measurable $E\subseteq N$ with $\nu(E)<\infty$  we obviously have
$$\int_E e^{\alpha|u|^{p'}}d\nu\le \int_{\{|u|\ge1\}} e^{\alpha|u|^{p'}}d\nu+e^\a\nu(E).$$
 From now on we will only focus on   exponential integrals over sets of finite measure,
since all the inequalities that involve regularized exponentials  stated in this paper can be deduced at once from this case, just by appealing to the lemma above.
In particular, under the  constraint $\|u\|_p\le 1$ it would be enough to prove uniform boundedness of the exponential integral over the set $E=\{|u|\ge 1\}$. This observation has been used before in connection with Moser-Trudinger inequalities in regularized form, see for example [LL2]. In other papers the same idea was used in different forms, after replacing $u$ by its radial decreasing rearrangement, in which case $E$ becomes a ball (see for example  [A], [Cao], [MS1], [Og],  [Ruf].)

To get to the heart of the matter, we now state and prove a  Lemma, which is   very elementary in nature, yet crucial  in the proof of   Theorem 5  and other results of this paper:\medskip

\proclaim Lemma 10. Let  $(N,\nu)$ be a measure space, and let $V,Z$ be  vector spaces of measurable functions (real, complex or vector valued) on a measurable set $E\subseteq N$. Let $L:V\to Z$ be an operator such that $L(\lambda f)=\lambda Lf$ for any $f\in V$ and any $\lambda\ge~0$,  and let $p:V\to [0,\infty]$ be a seminorm. Finally let $\b>1$ and $\b'=\ds{\b\over\b-1}$. \hfill\break 
\indent If there exists $c_0$ such that for a fixed subset $V_0\subseteq V$
$$\int_E \exp\bigg[{1\over A}|Lf(x)|^{\beta}\bigg]d\nu(x)\le c_0,\qquad \forall f\in V_0,\, p(f)\le 1$$
then, for each $\tau>0$\def\bbp{{\b\over\bp}} and for each $f\in V_0$ with $p(f)\le 1$ we have 
$$(A)\qquad\quad  \int_E \exp\bigg[{1\over A}\big(|L f(x)|+\tau\big)^{\beta}\bigg]d\nu(x)\le c_0 \exp\bigg[{1\over A} \bigg({\tau^{\bp}\over1-p(f)^\bp}\bigg)^{\bbp}\,\bigg]$$
and 
$$(B)\qquad \qquad \int_E \exp\bigg[{1\over A}\Big(|L f(x)|+\tau\big(1-p(f)^\bp\big)^{1\over\bp}\Big)^\b\,\bigg]d\nu(x)\le c_0e^{ \tau^\b/A}.$$
\par

\pf Proof. From H\"older's inequality we have 
$$a\theta^{\bpi}+b(1-\theta)^\bpi\le (a^\b+b^\b)^\bi,\qquad a,b\ge0,\,0\le\theta\le 1,\eqdef{20x}$$
from which estimate (B) follows, with $\theta=p(f)^\bp>0$, $a=L\big(f/ p(f)\big),\, b=\tau$. Clearly (A) is just another way of writing (B), since $\tau$ can be arbitrary.$\qquad /\!/\!/$ 
 
   \smallskip

As it is apparent from the proof, there is nothing peculiar about exponential integrability in this lemma, and (A) and (B) are clearly equivalent. Nonetheless, we find it convenient to have the estimates  in  (A) and (B) explicitly stated as above, since they will be used directly several times.
The first main application of the Lemma is in Theorem~5,  with $V=\{f\in L^\nsa(\Rn),\,\supp f $ compact$\}$, $Z=\{u:\R^n\to\R^m\,\,{\hbox{a.e. finite}}\}$, $\b=n/\a$, $p(f)=\|f\|_\na$, and $E\subseteq\Rn$ with $\nu=$Lebesgue measure (and $|E|\le 1$), and $L=T$.

For the benefit of the reader we will now summarize the main strategy behind the proof of Theorem 5. We start from the Adams inequality
$$\int_E\exp\bigg[{1\over A_g}\,|T f_1(x)|^{n\over n-\a}\bigg]dx\le C\eqdef{ADAMSf1}$$
valid for $|E|\le 1$ and all functions $f_1$ with $\|f_1\|_\na\le1$ with $|\supp f_1|\le \kappa$, for a fixed constant $\kappa$ (this follows from \eqref{12h}). The main idea is that the inequality in   Theorem 5  is true if we can  write $T f=T f_1+T f_2$, where $f=f_1+f_2$, $\,|\supp f_1|\le \kappa$, and where the additive perturbation $T f_2$ satisfies either
$$a)\;\;|T f_2(x)|\le C, \quad {\hbox {if}}\quad  x\in  E\qquad {\hbox{and}}\quad \|f_1\|_\na\le \theta_n<1$$
or
$$b)\;\;|T f_2(x)|\le C\big(1-\|f_1\|_\na^\na\big)^\asn, \quad {\hbox{if}}\quad  x\in E\qquad {\hbox{and}}\quad 0<\theta_n\le \|f_1\|_\na\le1, $$
 $\theta_n$ being a suitable explicit constant depending only on $n$. It is clear that in either case a) or b) one can apply (A) or (B) respectively, to derive the desired inequality.   The original set $E$ and the original function $f$  will  be suitably split so as to reduce matters to estimates a) and b). To this end, we will consider several scenarios  depending on where the $L^\na$ masses of $f,\, f^+,\,f^-$  are concentrated, and on where the potential is pointwise positive. An estimate for $T f_2$ as in a) will follow if the  function $f_2$ is either pointwise  small, or if its  support is ``well separated" from $E$; this will be a consequence of Lipschitz estimates for $T f$. An estimate as in b) will instead occur, roughly speaking,  when both $E$ and the mass of $f$ are concentrated in a fixed ball, and it will be the most critical case of the proof, the only one where Ruf's condition \eqref{11} is needed.

We now establish some regularity estimates for the operator $T$. 

\proclaim Lemma 11. Let $ f\in L^{\nsa}(\Rn)$, compactly supported and $F\subseteq \Rn$ a closed set such that  either
\smallskip\noin
(i) $|\supp f|\le 1$ and dist$ (F,  \supp f)\ge R\ge1$
\smallskip
\noin or
\smallskip\noin
(ii) $\,\supp f\subseteq B_{2R}^c$, $F\subseteq B_{R}$, 
\smallskip\noin
then $T f$ is Lipschitz on $F$, in particular there exists $D=D(n,\a,K)$ such that 
$$|T f(x_1)-T f(x_2)|\le{D\over R}\| f\|_\na |x_1-x_2|,\qquad x_1,x_2\in F.\eqdef{21}$$
Moreover, if $x_0\in F$ is so that $\max_{x\in F} |T  f(x)|=|T  f(x_0)|,$ then for $R\ge 2$
$$\max_{x\in F} |T  f(x)|\le \bigg(\avg_{B_1(x_0)}|T  f(x)|^\nsa dx\bigg)^{\asn}+{2D\over R}\| f\|_\na.\eqdef{22}$$\par
\pf Proof. For $x_1,x_2\in F$ and $y\in \supp f$, using the regularity estimate \eqref{rl3} we have 
$$\big|K(x_1-y)-K(x_2-y)\big|\le C |x_1-x_2|\cases {R^{\a-n-1} & in case (i)\cr
2^{n+1-\a}|y|^{\a-n-1} & in case (ii).\cr}$$
Hence,   in case (i) 
$$|T f(x_1)-T f(x_2)|\le C|x_1-x_2|R^{\a-n-1}\int_{\supp f}|f(y)|dy\le {D\over R}\|f\|_\na |x_1-x_2|,$$
whereas  in case (ii)
$$\eqalignno{&|T f(x_1)-T f(x_2)|\le C|x_1-x_2|2^{n+1-\a}\int_{B_{2R}^c}|f(y)||y|^{\a-n-1}dy & \eqdef{22k}\cr&\le C |x_1-x_2|
2^{n+1-\a}\|f\|_{\na}\bigg(\int_{B_{2R}^c}|y|^{-n-\nna}dy\bigg)^{{n-\a\over n}}\le{D\over R}\|f\|_\na |x_1-x_2|.\cr}$$
Note that using \eqref{rl2}
$$|K(x-y)|\le C|x-y|^{\a-n}\le C|x|^{\a-n},\qquad {\hbox { if }}x\neq0,\; |x|\ge {|y|\over2}$$
and since $\supp f$ is  compact
$|T f(x)|\le C_f |x|^{\a-n}$ for large $|x|$ (where $C_f$ depends on $f,\a,n, K$). Hence, the supremum of $|T f|$ over $F$ is attained in $F$ in either  case (i) or (ii). Let $M=\max_{x\in F} |T f(x)|=|T f(x_0)|$, some $x_0\in F$. After a moment's reflection the reader should realize that if $R\ge2$ then estimate \eqref{21} holds also for all $x_1,x_2$ in the set $F\cup B_1(x_0)$ (whose distance from $\supp f$ is at least $R-1$), by possibly enlarging slightly the constant~$D$.
 If $M\le {2D\over R}\|f\|_\na$ then \eqref{22} is true. If $M>{2D\over R}\|f\|_\na$ then  using \eqref{21} for  $x_1,x_2\in F\cup B_1(x_0)$ we get $|T f(x)|\ge M-{D\over R}\|f\|_\na\ge0$ for $x\in B_1(x_0)$ and 
$$\bigg( M-{D\over R}\|f\|_\na\bigg)^\nsa |B_1|\le \int_{B_1(x_0)} |T f(x)|^\nsa dx$$
which yields \eqref{22}.\endpf

\bigskip \proclaim Lemma 12 . If $f$ is compactly supported and $|f|\le 1$ on $\Rn$, then  there exists $D=D(n,\a)$ such that for $\alpha\neq 1$
$$|T f(x_1)-T f(x_2)|\le D(1+\|f\|_\na)\max\big\{|x_1-x_2|^\a,|x_1-x_2|\big\},\,\qquad x_1,x_2\in \Rn.\eqdef{23}$$
If $\alpha=1$ there exists $D=D(n)$ so that
$$|T_1 f(x_1)-T_1 f(x_2)|\le D(1+\|f\|_\na)|x_1-x_2|\Big(1+\log^+{1\over|x_1-x_2|}\Big),\,\qquad x_1,x_2\in \Rn,\;x_1\neq x_2.\eqdef{24}$$
\smallskip If $x_0\in\Rn$ is so that $\max_{x\in \sR^n} |T f(x)|=|T f(x_0)|$ then
$$\|T f\|_\infty\le \bigg(\avg_{B_1(x_0)}|T f(x)|^\nsa dx\bigg)^{\asn}+{2D}(1+\|f\|_\na).\eqdef{25}$$\par

 \par
\pf Proof. For each $x_1,x_2\in\Rn$ we have
$$\eqalign{&T f(x_1)-T f(x_2)= \!\!\!\!\mathop\int\limits_{|y-x_1|\le {1\over3}|x_1-x_2|}\hskip-2em K(x_1-y)f(y)dy+\hskip-1em \mathop\int\limits_{|y-x_2|\le {1\over3}|x_1-x_2|}\hskip-2em K(x_2-y)f(y)dy\cr&+ \hskip-1em\mathop\int\limits_{{|x_1-y|>|x_2-y|\atop |y-x_2|> {1\over3}|x_1-x_2|}}\hskip-2em\Big(K(x_2-y)-K(x_1-y)\Big)f(y)dy+\hskip-1em\mathop\int\limits_{{|x_2-y|>|x_1-y|\atop |y-x_1|> {1\over3}|x_1-x_2|}}\hskip-2em\Big(K(x_1-y)-K(x_2-y)\Big)f(y)dy.\cr}$$
then,  using \eqref{rl2}, \eqref{rl3}  we get
$$\eqalign{|T f(x_1)&-T f(x_2)|\le 2C\omega_{n-1}\int_0^{ {1\over3}|x_1-x_2|}r^{\a-1}dr+2L|x_1-x_2|\hskip-3em\mathop\int\limits_{{|y-x_1|<y-x_2|\atop  {1\over3}|x_1-x_2|<|y-x_1|\le 1}}\hskip-2em |x_1-y|^{\a-n-1}dy\cr&\hskip12em+2L|x_1-x_2|\mathop\int\limits_{{|y-x_1|\ge 1}} |x_1-y|^{\a-n-1}|f(y)|dy\cr&\le
{2C\omega_{n-1}\over \a 3^\a}|x_1-x_2|^\a+2\omega_{n-1}L|x_1-x_2|\int_{{1\over3}|x_1-x_2|}^1 r^{\a-2}dr\cr&\hskip12em +2L\Big({n-\a\over n}\Big)^{\asn}\|f\|_\na|x_1-x_2|\cr
}$$
and this proves \eqref{23} and \eqref{24}. Clearly $T f$ is continuous on $\Rn$ and $|T f(x)|\le C|x|^{a-n}$ for large $|x|$, so  $|T f(x)|$ has a maximum at some $x_0$. Estimate \eqref{25} is obtained as in Lemma~11.\endpf

\proclaim Lemma 13. If $f\in L^{\nsa}(\Rn)$, $\,|\supp f|<\infty
$, then there is $C=C(n,\a,K)$ such that  for all $E\subseteq \Rn$ with $|E|<\infty$
$$\bigg(\int_E |T f|^{\nsa}\bigg)^\asn\le C(|E|^{\asn}+|\supp f|^\asn) \|f\|_\na.\eqdef{25a}$$
\par\pf Proof. This  follows at once from O'Neil's inequality. If $f^*:(0,\infty)\to[0,\infty)$ denotes the nonincreasing rearrangement of a measurable function $f:\Rn\to \R$ and $f^{**}(t)=t^{-1}\int_0^t f^*(u)du$, then 
$$(T f)^{**}(t)\le t K^{**}(t)f^{**}(t)+\int_t^\infty K^*(t)f^*(t)dt\eqdef{oneil}$$
(See [ON], and also [FM1, Lemma 2] for an improvement). On the other hand, since $|K(x)|\le C|x|^{\a-n}$ we get $K^*(t)\le C  t^{\asn-1}$, and the same estimate holds for $K^{**}(t)$. Hence, for $t\le |E|$
$$(Tf)^{**}(t) \le C t^\asn f^{**}(t)+C\int_t^{|\supp f|} s^{\asn-1}f^*(s)ds\le C(|E|^\asn+|\supp f|^\asn) f^{**}(t),$$
if  for $t\le |\supp f|$, whereas if $t>|\supp f|$ the integral in \eqref{oneil} is 0 and the above  inequality is still valid, thereby proving the lemma.
\endpf

 Let now $f$ be compactly supported and such that 
$$\|f\|_\na^\na\le 1,\qquad \|T f\|_\na^\na\le1.\eqdef{25aa}$$
From now on we let

$$f=f_\ell+f_s,\qquad \quad f_\ell(x)=\cases{f(x) & if $|f(x)|\ge 1$\cr 0 & if $|f(x)|<1$.\cr}\eqdef{26}$$
Obviously   $\supp f_\ell$ is compact, $|f_\ell|\ge1$ in $\supp f_\ell$ and $|\supp f_\ell|\le1.$ Also, $|f_s|\le 1$ on $\Rn$, and $\supp f_s$ is compact, with no control on its measure.

\proclaim Lemma 14 . Under condition \eqref{25aa} we have
$$\|T f_s\|_{L^\infty(\sR^n)}\le C.\eqdef{27}$$
If $|F|<\infty$ then 
$$\|T\big( f_s\chi_F\big)\|_{L^\infty(\sR^n)}\le C(1+|F|^{\asn}).\eqdef{27a}$$
and if $|F^c|<\infty$ then 
$$\|T\big( f_s\chi_F\big)\|_{L^\infty(\sR^n)}\le C(1+|F^c|^{\asn}).\eqdef{27b}$$
\par
\pf Proof. From Lemma 12  we have
$$\|T f_s\|_\infty \le 2D(1+\|f_s\|_\na)+\bigg(\avg_{B_1(x_0)} |T f_s|^\nsa\bigg)^{\asn}$$
where $x_0$ is a  maximum for $|T f_s|$.
By Lemma 13 and \eqref{25aa}
$$\bigg(\int_{B_1(x_0)} |T f_s|^\nsa\bigg)^{\asn}\le \bigg(\int_{B_1(x_0)} |T f|^\nsa\bigg)^{\asn}+\bigg(\int_{B_1(x_0)} |T f_\ell|^\nsa\bigg)^{\asn}\le 1+C\|f_\ell\|_\na\le C.$$
If $|F|<\infty$ then \eqref{27a} follows from \eqref{25}, \eqref{25a}, \eqref{25aa}. If $|F^c|<\infty$ then
just write $f_s\chi_F=f_s-f_s\chi_{F^c}$.
\endpf

\bigskip \centerline{\bf 4. Proof of   the inequalities  in  Theorem 5 }\bigskip

It is enough to prove \eqref{12} in  the case $|E|\le 1$, since \eqref{11} implies $|\{|Tf|\ge1\}|\le 1$, and we can use Lemma 9.
Let then   $E\subseteq \R^n$ be measurable and with $|E|\le 1,$ and let $f$ be compactly supported and satisfying \eqref{25aa}. 

After  suitable translations we can  assume that
$$\int_{x_i\le0} |f(x)|^\nsa dx=\int_{x_i\ge0}|f(x)|^\nsa dx={1\over2}\; \nf,\qquad i=1,2,...,n.\eqdef{29}$$
Define the $2n$ half-spaces
$$H_i^+=\{x\in\Rn: \;x_i\ge4\},\qquad H_i^-=\{x\in\Rn: \; x_i \le -4\}\eqdef{30}$$
\medskip
We organize the proof in 6  cases:
\medskip
\noin{\bf \underbar{Case 1}:} There exists $i\in\{1,2...,n\}$ and  a half-space $H\in \{H_i^+,H_i^-\}$ such that 
$$\int_H |f(x)|^\nsa dx\ge {1\over4n}\nf.\eqdef{31}$$
\medskip
\noin{\bf \underbar{Case 2}:} 
$$\int_{B_{4\sqrt n}}|f(x)|^\nsa dx\ge{1\over2}\nf,\quad{\hbox {and}}\quad \|f_\ell\|_\na^\na\le{3\over4}\nf.\eqdef{32}$$
\medskip
\noin{\bf \underbar{Case $\bf {2^+}$}:}  The kernel $K$ is {\it {nonnegative}} and 
$$\int_{B_{4\sqrt n}}|f(x)|^\nsa dx\ge{1\over2}\nf,\quad \|f_\ell^+\|_\na^\na\le{3\over4}\nf,\;\;\|f_\ell^-\|_\na^\na\le{3\over4}\nf.\eqdef{32a}$$

\medskip
\noin{\bf \underbar{Case 3}:}     $E\subseteq B_{16\sqrt n}^c$ and 
$$\int_{B_{4\sqrt n}}|f(x)|^\nsa dx\ge{1\over2}\nf,\quad \|f_\ell\|_\na^\na\ge{3\over4}\nf.\eqdef{32aa}$$


\medskip\noin{\bf \underbar{Case $\bf 4$}:}  $E\subseteq B_{16\sqrt n}$ and conditions \eqref{11}, \eqref{12q} hold.\medskip


\medskip\noin{\bf \underbar{Case $\bf 4^+$}:}   $E\subseteq B_{16\sqrt n}$, the kernel $K$ is {\it {nonnegative}} and 
$$\int_{B_{4\sqrt n}}|f(x)|^\nsa dx\ge{1\over2}\nf,\quad \|f_\ell^+\|_\na^\na\ge{3\over4}\nf\;\;{\hbox{ OR}}  \;\;\|f_\ell^-\|_\na^\na\ge{3\over4}\nf.\eqdef{32b}$$

Within this case we will consider the  following subcases:\medskip
\smallskip\item{\bf{(i)}} $T f(x)\le0$ for any $x\in E$;\smallskip
\item{\bf{(ii)}}  $T(f\chi_{B_{32\sqrt n}}^{})(x)$ and $T(f\chi_{B_{32\sqrt n}^c}^{})(x)$ have opposite sign, and $Tf(x)\ge0$, for any $x\in E$;\smallskip
\item{\bf{(iii)}}  $T(f\chi_{B_{32\sqrt n}}^{})(x)\ge0$ and $T(f\chi_{B_{32\sqrt n}^c}^{})(x)\ge0$ for any $x\in E$, and condition \eqref{11} holds.\smallskip

\medskip
It is easy to see that once    Theorem 5  is proved in the above cases then it is proved in full generality. Indeed,  if Case 1 is not verified then $\int_{B_{4\sqrt n}}|f|^{\nsa}\ge \half\nf$, so that from Cases 1,2,3 the theorem follows when $E\subseteq B_{16\sqrt n}^c$. For arbitrary $E$ write $E=(E\cap  B_{16\sqrt n})\cup(E\cap  B_{16\sqrt n}^c)$ and the theorem follows if \eqref{12q} is assumed. Similarly the theorem follows from 1, $2^+$, $3$, $4^+$ if $K$ is nonnegative, and the case $K$ nonpositive obviously follows as well.

It is worth emphasizing  that  {\it cases 1, 2, 3 and 4  hold for vector-valued kernels with arbitrary sign}, and that {\it the stronger ``Ruf condition" \eqref{11} is only used  in cases $4$ and $4^+$ }(iii).

We will now prove the main estimate \eqref{12}  in these cases.
\medskip
\noin\underbar {\sl Proof of \eqref{12} in Case 1.}  Suppose  WLOG that 
$$\int_{x_1\ge4} |f(x)|^\nsa dx\ge{1\over4n}\nf\eqdef{33}$$
and write
$$\int_E\exp\bigg[{1\over A_g}\,|T f(x)|^{n\over n-\a}\bigg]dx=\int_{E\cap\{x_1\le 2\}}+\int_{E\cap\{x_1\ge2\}}=I+II.\eqdef {34}$$
To estimate $I$ write
$$T f=T\big(f_\ell\chi_{\{x_1\le 4\}}^{}\big)+T\big(f_\ell\chi_{\{x_1\ge 4\}}^{}\big)+T f_s.\eqdef{35}$$
From Lemmas 11 and 13 we get 
$$\big|T\big(f_\ell\chi_{\{x_1\ge 4\}}^{}\big)(x)\big|\le \bigg(\avg_{B_1(x_0)}\big|T\big(f_\ell\chi_{\{x_1\ge 4\}}^{}\big) (x)\big|^\nsa dx\bigg)^{\asn}+2D\|f_\ell\|_\na\le C,\quad  
\forall x\in\{x_1\le 2\}$$
(here $x_0$ is a maximum for $\big|T\big(f_\ell\chi_{\{x_1\ge 4\}}^{}\big) (x)\big|$ in $\{x_1\le 2\}$.)
From Lemma 14  we then have
$$|T f(x)|\le \big|T\big(f_\ell\chi_{\{x_1\le 4\}}^{}\big) (x)\big|+\tau,\qquad \forall x\in\{x_1\le 2\}$$
some $\tau$ depending only on $\a,n,K$, and (A) of Lemma 10 and \eqref{ADAMSf1} (with $f_1=f_\ell\chi_{\{x_1\le 4\}}^{}$)  imply that $I\le C$, since $$\|f_\ell\chi_{\{x_1\le 4\}}^{}\|_\na^\na\le\|f\chi_{\{x_1\le 4\}}^{}\|_\na^\na\le \Big(1-{1\over 4n}\Big)\nf<1-{1\over 4n}<1.$$
The estimate of $II$ is similar, this time write
$$T f=T\big(f_\ell\chi_{\{x_1\le 0\}}^{}\big)+T\big(f_\ell\chi_{\{x_1\ge 0\}}^{}\big)+T f_s.\eqdef{36}$$
and Lemmas 11,  13, 14, imply
$$|T f(x)|\le \big|T\big(f_\ell\chi_{\{x_1\ge 0\}}^{}\big) (x)\big|+\tau,\qquad \forall x\in\{x_1\ge 2\}$$
some $\tau$ depending only on $\a,n,K$, and (A) of Lemma 10  and \eqref{ADAMSf1} imply that $II\le C$, since 
$$\|f_\ell\chi_{\{x_1\ge0\}}^{}\|_\na^\na\le\|f\chi_{\{x_1\ge 0\}}^{}\|_\na^\na= {1\over2}\nf<{1\over2}.$$
\medskip

\noin\underbar {\sl Proof of \eqref{12} in Case 2.} \smallskip
Assume \eqref{32} and write
$$\int_E\exp\bigg[{1\over A_g}|T f(x)|^\nna\bigg] dx\le\int_{E}\exp\bigg[{1\over A_g}\Big(T_{|K|} |f_\ell|(x)+|T f_s(x)|\Big)^\nna \bigg]dx$$
with $T_{|K|}f=|K|*f$,  and \eqref{12} follows from \eqref{27a},   \eqref{ADAMSf1} and (A) of Lemma 10 applied to $f_\ell$, since $\|f_\ell\|_\na^\na\le {3\over4}$.
\bigskip
\noin\underbar {\sl Proof of \eqref{12} in Case $2^+$.} \smallskip
Assume \eqref{32a} and  write
$$\eqalign{\int_E\exp\bigg[{1\over A_g}|&T f(x)|^\nna\bigg] dx\le\int_{E^+}\exp\bigg[{1\over A_g}\Big(T f_\ell^+(x)+|T f_s(x)|\Big)^\nna \bigg]dx+\cr&+\int_{E^-}\exp\bigg[{1\over A_g}\Big(T f_\ell^-(x)+|T f_s(x)|\Big)^\nna\bigg] dx\cr}$$
where $E^{\pm}=\{x\in E:\,T f_\ell^\pm(x)\ge T f_\ell^\mp(x)\}$, and \eqref{12} follows from (A) and  \eqref{ADAMSf1} , applied to $f_\ell^\pm$.

\bigskip
\noin\underbar {\sl Proof of \eqref{12} in Case 3.} Suppose $E\subseteq  B_{16\sqrt n}^c$ and that the estimates in  \eqref{32aa} hold.  
For $x\in E$ we then have 
$$|T f(x)|\le \big|T\big(f_\ell\chi_{B_{8\sqrt n}}^{}\big)(x)\big|+\big
|T\big(f_\ell\chi_{B_{8\sqrt n}^c}^{}\big)(x)\big|+|T f_s(x)|\le\big|T\big(f_\ell\chi_{B_{8\sqrt n}^c}^{}\big)(x)\big|+ C$$
from Lemmas 11-14.
  On the other hand, since $\big\|f\chi_{B_{4\sqrt n}}^{}\big\|_\na^\na\ge{1\over2}\nf$, it must be that $\big\|f_\ell\chi_{B_{8\sqrt n}^c}^{}\big\|_\na^\na\le{1\over2}\nf<{1\over2}$, and $\big|\supp\big(f_\ell\chi_{B_{8\sqrt n}^c}^{}\big)\big|\le 1$, so that our estimate \eqref{12}  follows again from (A) of Lemma 10 and  \eqref{ADAMSf1},  applied  to $f_\ell\chi_{B_{8\sqrt n}^c}^{}$.
\bigskip
\noin\underbar {\sl Proof of \eqref{12} in Case $4$.} \smallskip
Let  $E\subseteq  B_{16\sqrt n}$ and suppose that $f$ satisfies both  \eqref{11}  and \eqref{12q}. Estimate \eqref{12q} implies that 
$$|T f(x)|\le \int_{B_{32\sqrt n}}|K(x-y)|\,|f(y)|dy+C_1\|T f\|_\na,\qquad  x\in B_{16\sqrt n}.$$
(if $x\in B_{16\sqrt n}$, pick $a\in B_{16\sqrt n}$ such that $|x-a|\le 1$, use \eqref{12q}, and enlarge the domain of integration.)
The Ruf condition $\nf+\ntf\le1$  gives, with  $T_{|K|}f=|K|*f$,
$$ \eqalign{|Tf(x)|&\le  T_{|K|}\big(|f|\chi_{B_{32\sqrt n}}^{}\big)+C(1-\|f\|_\na^\na)^\an\cr& \le  T_{|K|}\big(|f|\chi_{B_{32\sqrt n}}^{}\big)+C(1-\|f\chi_{B_{32\sqrt n}}^{}\|_\na^\na)^\an\cr}\qquad x\in B_{16\sqrt n,}\eqdef{usaruf1}$$
and  Lemma 10 (B) yields inequality \eqref{12}.
\bigskip
\noin\underbar {\sl Proof of \eqref{12} in Case $4^+$.} \smallskip
It is clear that   in   \eqref{32b} it is  enough to assume that  $\|f_\ell^+\|_\na^\na\ge{3\over4}\nf$, in which case we have
$$\|f_\ell^-\|_\na^\na\le {1\over4}\nf,\eqdef{39}$$
and
$$\int_{B_{4\sqrt{n}}} (f_\ell^+)^{\nsa} \ge\irn|f_\ell^+|^\nsa-\int_{B_{4\sqrt n}^c}|f|^{\nsa}\ge{1\over4}\nf.\eqdef{40}$$

In case (i) we have $Tf(x)\le 0$ in $E$, hence, since $g\ge0$, 
$$\eqalign{\int_E\exp\bigg[{1\over A_g}|&T f(x)|^\nna\bigg] dx\le\int_{E}\exp\bigg[{1\over A_g}(T f_\ell^-(x)+|T f_s(x)|)^\nna\bigg] dx\cr}\le C,$$
using again (A)  and \eqref{ADAMSf1} applied to $f_\ell^-$. 

In case (ii) write
$$\int_E\exp\bigg[{1\over A_g}|T f(x)|^\nna\bigg] dx=\int_{E}\exp\bigg[{1\over A_g}
\Big(T\big(f\chi_{B_{32\sqrt n}}^{}\big)+T\big(f\chi_{B_{32\sqrt n}^c}^{}\big)\Big)^\nna\bigg]dx.$$

If $T\big(f\chi_{B_{32\sqrt n}}^{}\big)\ge0 $ and $T\big(f\chi_{B_{32\sqrt n}^c}^{}\big)\le0$ then 
$$\int_E\exp\bigg[{1\over A_g}|T f(x)|^\nna\bigg] dx\le\int_E\exp\bigg[{1\over A_g}\Big(T\big(f\chi_{B_{32\sqrt n}}^{}\big)(x)\Big)^\nna \bigg]dx\le C$$
by the original Adams inequality \eqref{12h}.

If instead $T\big(f\chi_{B_{32\sqrt n}}^{}\big)\le0 $ and $T\big(f\chi_{B_{32\sqrt n}^c}^{}\big)\ge0$ then 
$$\int_E\exp\bigg[{1\over A_g}|T f(x)|^\nna\bigg] dx\le\int_E\exp\bigg[{1\over A_g}\Big(T\big(f\chi_{B_{32\sqrt n}^c}^{}\big)(x)\Big)^\nna \bigg]dx\le C$$
since on $E\subseteq B_{16\sqrt n}^{}$ we have
$$0\le T\big(f\chi_{B_{32\sqrt n}^c}^{}\big)(x)\le |T\big(f_\ell\chi_{B_{32\sqrt n}^c}^{}\big)(x)|+
|T\big(f_s\chi_{B_{32\sqrt n}^c}^{}\big)(x)|\le C$$
from Lemmas 11-14.

In case (iii), the most critical situation, let us assume the Ruf condition
$\nf+\ntf\le1$
 and write 
$$\eqalign{\int_E\exp\bigg[{1\over A_g}|T f(x)|^\nna\bigg] dx\le\int_{E}\exp\bigg[&{1\over A_g}\bigg(T\big(f^+\chi_{B_{32\sqrt n}}^{}\big)(x)+\cr&+T\Big((f_\ell^++f_s)\chi_{B_{32\sqrt n}^c}^{}\Big)(x)\bigg)^\nna \bigg]dx.\cr}\eqdef{11h}$$

By Lemma 11 , for $x\in E\subseteq B_{16\sqrt n}^{}$
$$\eqalign{0\le T\Big((f_\ell^++f_s)\chi_{B_{32\sqrt n}^c}^{}\Big)(x)\le\bigg(\avg_{\!\!\!B_1(x_0)} \Big|T\Big((f_\ell^++f_s)\chi_{B_{32\sqrt n}^c}^{}&\Big)\Big|^{\nsa}dx\bigg)^{\asn}+\cr&+C\big\|(f_\ell^++f_s)\chi_{B_{32\sqrt n}^c}^{}\big\|_\na\cr}$$
where $x_0$ is a maximum point for $T\Big((f_\ell^++f_s)\chi_{B_{32\sqrt n}^c}^{}\Big)$ on 
$\overline{B}_{16\sqrt n}$. 

We have 
$$\eqalign{\bigg(\int_{B_1(x_0)} \Big|T\Big((f_\ell^+&+f_s)\chi_{B_{32\sqrt n}^c}^{}\Big)\Big|^{\nsa}dx\bigg)^{\asn}\le\bigg(\int_{B_1(x_0)} \Big(T\big(f_\ell^++f_s^++f_s^-\chi_{B_{32\sqrt n}^c}^{}\big)\Big)^{\nsa}dx\bigg)^{\asn}\cr&\le \bigg(\int_{B_1(x_0)} |T f|^\nsa dx\bigg)^\asn+\bigg(\int_{B_1(x_0)} \Big|T\big(f_\ell^-+f_s^-\chi_{B_{32\sqrt n}}^{}\big)\Big|^{\nsa}dx\bigg)^\asn\cr&\le\|T f\|_\na+C\Big(\|f_\ell^-\|_\na+\big\|f_s^-\chi_{B_{32\sqrt n}}^{}\big\|_\na\Big)\cr}$$
(by Lemma 13 ). Hence, for $x\in E$ 
$$\eqalign{0&\le T\Big((f_\ell^++f_s)\chi_{B_{32\sqrt n}^c}^{}\Big)(x)\le\cr& \le C\Big(\|T f\|_\na+\|f_\ell^-\|_\na+\big\|f_s^-\chi_{B_{32\sqrt n}}^{}\big\|_\na+\big\|f_\ell^+\chi_{B_{32\sqrt n}^c}^{}\big\|_\na+\big\|f_s\chi_{B_{32\sqrt n}^c}^{}\big\|_\na\Big)\cr&\le C \Big(\|T f\|_\na^\na+\|f_\ell^-\|_\na^\na+\big\|f_s^-\chi_{B_{32\sqrt n}}^{}\big\|_\na^\na+\big\|f_\ell^+\chi_{B_{32\sqrt n}^c}^{}\big\|_\na^\na+\big\|f_s\chi_{B_{32\sqrt n}^c}^{}\big\|_\na^\na\Big)^\an\cr&\le C\Big(1-\big\|f_\ell^+\chi_{B_{32\sqrt n}}^{}\big\|_\na^\na-\big\|f_s^+\chi_{B_{32\sqrt n}}^{}\big\|_\na^\na\Big)^\an= C\Big(1-\big\|f^+\chi_{B_{32\sqrt n}}^{}\big\|_\na^\na\Big)^\an\cr}\eqdef{usaruf2}$$
where in the second to  last estimate we used that $\ntf\le1-\nf$.

By applying (B) of Lemma 10 to the function $f^+\chi_{B_{32\sqrt n}}^{}$ we obtain that the integral in \eqref{11h} is bounded by  a constant $C$. This concludes the proof of the inequality part of Theorem~5. 

The sharpness of the exponential constant will be proved in Section  6.

\endpf

\noin{\bf Note.} We would like to emphasize that the precise steps where Ruf's condition is used are in \eqref{usaruf1} and \eqref{usaruf2}. Those steps also make it clear that the reason why the proof of   Theorem 5  fails if one uses the condition $\|f\|_\na^{q\na}+\|Tf\|_\na^{q\na}\le1$ is that the inequality $\|f\|_\na^{qn/\a}\ge\|f\chi_B^{}\|_\na^{qn/\a}+\|f\chi_{B^c}^{}\|_\na^{qn/\a}$ is only true for $q\le1$ (being trivially an equality when $q=1$).

\bigskip
  
\centerline{\bf 5.   Proof of   the inequalities in  Theorem 6 }\bigskip

Assume the Adachi-Tanaka type condition 
$$\max\big\{\|f\|_\na,\|Tf\|_\na\big\}\le1$$
and let us prove \eqref{102}, i.e. for $0<\theta<1$
$$\int_E \exp\Big[{\theta\over A_g}|T f(x)|^\nna\Big]dx\le C_\theta(1+|E|).\eqdef{AT1}$$
This proof does not require the full force of   Theorem 5, but uses instead only (A) in Lemma~10 and Lemma 14. Indeed, we have $|T (\theta^{n-\a\over n} f)|\le |T (\theta^{n-\a\over n} f_\ell)|+C$ by Lemma 14, hence \eqref{AT1} follows at once  from (A) of Lemma 10 and \eqref{ADAMSf1}, since $1-\|\theta^\nna f_\ell\|_\na^\na\ge 1-\theta^{n-\a\over\a}$.
Alternatively one can just use
$$|Tf|^p\le\big||Tf_\ell|+C\big|^p\le (1+\e)|Tf_\ell|^p+C_\e$$
with $\theta(1+\e)<1$.
 Obviously  \eqref{AT1}  follows under the more restrictive  condition $$\|f\|_\na^{q\na}+\|T f\|_\na^{q\na}\le~1,\qquad q\ge 1,$$ but in the next section we show that under this condition with $1<q\le+\infty$ the inequality fails if $\theta=1$. 

In the case $K$ homogeneous we can proceed with a dilation argument, using Theorem~5. 
Assume the Adachi-Tanaka condition \eqref{103}. Clearly,  it suffices to prove  estimate \eqref{102} for $0<\theta_0\le\theta<1$.
For any $\lambda>0$ if $f_\lambda(x)=\lambda^\a f(\lambda x)$ then using \eqref{D1} 
$$\int_E \exp\Big[{\theta\over A_g}|T f(x)|^\nna\Big]dx=\lambda^n\int_{E/\lambda} \exp\Big[{1\over A_g}|T (\theta^{n-\a\over n} f_\lambda)(x)|^\nna\Big]dx\eqdef{D3}$$
and, using H\"older's inequality
$$\eqalign{\|\theta^{n-\a\over n} f_\lambda\|_\na^\na+\|T(\theta^{n-\a\over n} f_\lambda)\|_\na^\na&=\theta^{n-\a\over\a}\Big(\|f\|_\na^\na+\lambda^{-n}\|T f\|_\na^\na\Big)\cr&\le\theta^{n-\a\over\a}\big(1+\lambda^{-nq'}\big)^{1/q'}=1\cr}\eqdef{D4}$$
 for $\lambda=\lambda(\theta)=\big(\theta^{-q'{n-\a\over \a}}-1\big)^{-{1\over nq'}}\ge \lambda(\theta_0)=1$, if we choose $\theta_0=2^{-{n\over q'(n-\a)}}$.
We can then apply Theorem 5  to estimate \eqref{D3} with $C\lambda^n\le C(1-\theta)^{-1/q'}$.

\bigskip

\centerline{\bf 6. Proof of the sharpness statements in Theorems 5 and 6}
\bigskip

Now we show that the constant $ A_g^{-1}$ is best possible in \eqref{12}, and \eqref{102}, in the sense that if $|E|>0$ then  we can find a family of functions $\psi_\e\in L^{\na}(\Rn)$ such that $\|\psi_\e\|_\na^{q\na}+\|T \psi_\e\|_\na^{q\na}\le 1$, compactly supported, and 
$$\lim_{\e\to0^+}\int_{E\cap B_{\e/2}}\exp\bigg[{\theta\over A_g}\,|T \psi_\epsilon(x)|^{n\over n-\a}\bigg]dx=+\infty\eqdef {12x}$$
if either $q=1$, $\theta>1$ or $q>1, \theta=1$. For slightly more clarity we will  prove the result in the scalar case first, but the modifications in the vector-valued case are simple, and will be indicated after the proof of the scalar case. 

First notice that for a.e. $x\in E$ we have $|E\cap B_\e(x)|/|B_\e(x)|\to1$, as $\e\to0$, therefore we can assume WLOG that for some $\e_0>0$
$$|E\cap B_\e|\ge{\half} |B_1| \e^n,\qquad 0<\e<\e_0.\eqdef{12w}$$  
It is no big surprise that even in the case of the Riesz potential $K=I_\a$  a family of functions satisfying \eqref{12x}, with $A_g=|B_1|$,  will be obtained by a  suitable modification of the usual extremal Adams family
$$\phi_\e(y)=\cases{|y|^{-\a} & if $\e< |y|\le 1$\cr 0 & otherwise\cr},\qquad 0\le\e<1$$
It is clear that some modification is necessary for some values of $\alpha$,  due to the integrability requirements on $(I_\a \phi_\e)^\na$ at infinity. Indeed,  as $|x-y|^{\a-n}\sim |x|^{\a-n}$ when  $|x|$ is large, and this implies that for any $f$ compactly supported in $\Rn$,  $I_\a f$ cannot be in $L^\na(\R^n)$ for $n/2\le \alpha<n$,  if $\int_{\sR^n} f\neq0$. 

 Assume that $K$ is $n-$regular and $|y|\le \half |x|$, $x\neq0$.  From  Taylor's formula centered at $x$, for any positive integer $m\le n$, there exists $\theta=\theta(x,y,m)\in(0,1)$ such that 
$$K(x-y)=\sum_{j=0}^{m-1} p_j(x,y)+ R_m(x,y)\eqdef{g}$$
where 
$$\cases{p_j(x,y)=\ds{1\over j!}d^jK(x,-y) &  $j=0,1,...,m-1$
\cr R_m(x,y)=\ds{1\over m!}d^m K(x-\theta y,-y)\cr}\eqdef{Taylor1}$$
 and where $d^j K(x;y)$ denotes the $j-$th differential of $K$ at $x$ with increment $y$, i.e.
$$d^jK(x;y)=\mathop\sum\limits_{j_1+...+j_n=j}{j!\over j_1!...j_n!}\,{\p^jK\over \p x_1^{j_1}...\p x_n^{j_n}}(x)y_1^{j_1}...y_n^{j_n}.\eqdef{Taylor2}$$
Clearly,  $p_j(x,y)$ is a polynomial of order $j$ in the $y$ variable for $j\le m-1$, and for  $x\neq 0$
$$\cases{|p_j(x,y)|\le C|x|^{\a-n-j}|y|^j & for $j\le m-1,\,y\in\Rn,$ \cr
\cr |R_m(x,y)| \le C|x|^{\a-n-m}|y|^m  & for $|y|\le \half|x|, $\cr}\eqdef{stima}$$
for some constant $C$ independent of $y$ and $x$ in the given range.

\proclaim  Lemma 15. If $K$ is an $n-$regular Riesz-like kernel and $f\in L_c^1(\Rn)$, with $\supp f\subseteq B_r$. Then we have that $Tf\in L^\na(|x|\ge 2r)$ for $0<\a<{n\over2}$. If ${n\over2}\le \a<{n+m\over2}$,  $\,m=1,2...n$, then $Tf\in L^\na(|x|\ge2r)$ if $f$ is orthogonal to the family of homogeneous polynomials $\{p_j(x,\cdot), \,j=0,1,...,m-1,\;x\in \Rn\},$ in particular if $f$ has vanishing moments up to order $m-1$.
In any of the above situations we have
$$\int_{|x|\ge 2r}|Tf(x)|^\na dx\le C r^{2n-{n^2\over\a}}\|f\|_1^\na.\eqdef{normaT}$$\par

\noin{\bf Remarks.} 

\noin 1) As a consequence of the above Lemma and  Lemma 13 we have that if $f\in L_c^\na(\Rn)$ has enough vanishing moments, as specified above, then $Tf\in L^\na(\Rn)$. 

\noin{2)}  In the case  $K(x)=|x|^{\a-n}$ it is possible to evaluate $p_j$ explicitly via the binomial expansion:
$$p_j(x,y)=|x|^{\a-n-j}|y|^{-j}\mathop\sum\limits_{j/2\le k\le j}{{\a-n\over2}\choose {k}} {{k}\choose{2k-j}}(-2)^{ 2k-j}(x^*\cdot y^*)^{2k-j}.$$
With the aid of such formula it is easy  to check that the vanishing of all moments up to order $0,1,2$ is equivalent to the orthogonality to the functions $p_j(x,\cdot)$, $\,j=0,1,2$, as $x\in \Rn$; this fact is very  likely true for all  moments up to any given order.
\bigskip
\pf Proof of Lemma 15. Suppose that $\supp f\subseteq B_r$.  If $\a<n/2$ we use the estimate $|K(x-y)|\le C|x-y|^{\a-n}$ to deduce that $|Tf(x)|\le C|x|^{\a-n}\int_{|y|\le r}|f|$, for $|x|\ge 2r$, which implies \eqref{normaT}.

Let now  ${n\over2}\le \a<{n+m\over2}$ and let $f$ be  orthogonal to all of the $p_j(x,y)$ up to $j=m-1$. If $|x|\ge 2r$ then using \eqref{g},\eqref{stima}
 $$|Tf(x)|=\bigg|\int_{B_r} R_m(x,y)f(y)dy\bigg|\le C|x|^{\a-n-m}\int_{|y|\le r}|y|^m |f(y)|dy\le Cr^m |x|^{\a-n-m}\|f\|_1$$
and this implies \eqref{normaT}.\endpf

Next, let $\P_m$ be  the space of polynomials of degree up to $m$ in the unit ball $B_1$ of $\Rn$, a subspace of $L^2(B_1)$. Let $\{v_1,....,v_N\}$ be an orthonormal basis of $\P_m$, with $v_1=|B_1|^{-1/2}$. If $P_m$ denotes the projection of $L^2(B_1)$ onto $\P_m$, then $P_m$ has integral kernel
$$P_m(y,z)=\chi_{B_1}^{}(y)\sum_{k=1}^N v_k(y)v_k(z)\eqdef{proj}$$
which is pointwise bounded on $B_1\times B_1$ (with bound depending on $m$). The projection can be then extended to all of $L^1(B_1)$, and the function $f-P_m f$ is  orthogonal to all polynomials on $B_1$ of order up to $m$.

For $r>0$ the projection on the space of polynomials of degree up to $m$ on the ball $B_r=B(0,r)$ is given by the integral kernel $$P_m^r(y,z)=r^{-n}P_m\Big({y\over r},{z\over r}\Big).$$

\def\phier{\phi_{\e,r}}
For $0\le\e r<1\le r$ let $$\phier(y)=\cases{K(-y)|K(-y)|^{{\a\over n-\a}-1}  & if $\e r< |y|\le r$\cr 0 & otherwise\cr}.\eqdef{phie}$$
and  consider the functions on $B_r$
$$\wtilde \phier=\phier-P_{n-1}^r\phier,\eqdef{12g}$$
which are orthogonal to every polynomial of order up to (and including) $n-1$, on the ball $B_r$.

Note that \eqref{rl2}  implies $\|\phier\|_1\le C r^{n-\a}$, so that   $|P_{n}^r\phier|\le Cr^{-n}\|\phier\|_1\le Cr^{-\a}$ and
$$|\wtilde \phier(-y)|\le |K(y)|^{\a\over n-\a}+Cr^{-\a},\qquad \e r< |y|\le r.$$

Using the inequality $(a+b)^p\le a^p+p2^{p-1}(a^{p-1}b+b^p)$ we have
$$|\wtilde \phier(-y)|^{\na}\le |K(y)|^{\nna}+Cr^{-\a}|K(y)|+Cr^{-n},\qquad \e r< |y|\le r.$$
When $\e r< |y|\le 1$ we  use  \eqref{rl1} to obtain 
$$\eqalign{|\wtilde \phier(-y)|^{\na}&\le |K(y)|^{\nna}+C|K(y)|+C\cr&\le 
|y|^{-n}\big(|g(y^*)|+C|y|^\d\big)^\nna+C|y|^{\a-n}\big(|g(y^*)|+C|y|^\d\big)+C
\cr&
\le|g(y^*)|^\nna|y|^{-n}+C\big(|y|^{-n+\d}+|y|^{-n+{n\d\over n-\a}}+|y|^{\a-n}+1\big)\cr&=|g(y^*)|^\nna|y|^{-n}+\Psi_0(y)\cr}$$
where $\Psi_0\in L^1(|y|\le 1)$, and consequently
$$\eqalign{\|\wtilde \phier\|_\na^\na&\le\int_{\e r< |y|\le 1}\Big(|g(y^*)|^\nna|y|^{-n}+\Psi_0(y)\Big)dy+\int_{1\le |y|\le r}|K(y)|^\nna dy\cr&+C\int_{1\le |y|\le r}\big(r^{-\a}|y|^{\a-n}+r^{-n}\big)dy
=
 A_g\log{1\over( \e r)^n}+b_r+C\cr}\eqdef{12gg}$$
where 
$$b_r:=\int_{1\le |y|\le r}|K(y)|^\nna dy,\qquad b_r\le C\log r.$$

For later use we also note the following estimate from below:
$$ \int_{\e r< |y|\le r}|K(y)|^\nna dy= \int_{\e r<|y|\le 1}|K(y)|^\nna dy+ b_r\ge A_g\log{1\over(\e r)^n}-C+b_r\eqdef{r4}$$
which follows from
$$|K(y)|^\nna\ge \big||K_0(y)|-|K(y)-K_0(y)|\big|^\nna,$$
where $K_0(y)=g(y^*)|y|^{\a-n}$. Indeed, using \eqref{rl1} and  the elementary inequality  $|a-b|^p\ge |a|^p-p|a|^{p-1}|b|$, valid for all $a,b\in\R$ and $p\ge 1$  we obtain
 $$|K(y)|^\nna\ge|g(y^*)|^\nna|y|^{-n}-\Phi_0(y),\qquad \e r<|y|\le 1$$
where $\Phi_0(y)\in L^1(|y|\le 1)$.

Now let us  estimate $\|T\wtilde\phier\|_\na$.  
If $|x|\ge 2r$, then, since $\wtilde \phier$ is orthogonal to every $p_j$, $j\le n-1$, from Lemma 15 and  \eqref{normaT}  we have
$$\int_{|x|\ge 2r}|T\wtilde \phier(x)|^{\na}\le Cr^{2n-{n^2\over\a}}\|\wtilde \phier\|_1^\na\le Cr^{2n-{n^2\over\a}}r^{{(n-\a)}\nsa}=Cr^n.\eqdef{12t}$$
If  $|x|\le 2r$ we have
$$|TP_n^r\phier(x)|\le Cr^{-\a}\int_{|z-x|\le r}|z|^{\a-n}dz\le Cr^{-\a}\int_{|z|\le 4r}|z|^{\a-n}dz= C\eqdef{TP}$$
 $$|T\wtilde \phier(x)|\le C\int_{|y|\le r} |x-y|^{\a-n}|y|^{-\a}+C=C\int_{|y|\le 1} \Big|{x\over r}-y\Big|^{\a-n}|y|^{-\a}+C:=CI_\a\phi_{0}\Big({x\over r}\Big)+C.$$
But $\phi_{0}(y)=|y|^{-\a}\chi_{|y|\le 1}^{}$ is in  $L^p(\Rn)$ if $1<p<\nsa$, hence $I_\a\phi_{0}\in L^q(\Rn)$
with $q^{-1}=p^{-1}-\a/n$. If we pick any $p$ so that ${n\over2\a}<p<\nsa$ then $q>\nsa$ and so  $I_\a\phi_{0}\in L^{\na}(|x|\le 2)$. This shows 
$$\int_{|x|\le2r} |T\wtilde \phier(x)|^\na\le C \int_{|x|\le 2r} \Big|I_\a\phi_{0}\Big({x\over r}\Big)+C\Big|^\na dx= Cr^n$$
and therefore $$\|T\wtilde\phier\|_\na^\na\le C_1 r^n\eqdef{T}$$ with $C_1$ independent of $\e,r$.
Let us denote temporarily \def\ber{b_{\e,r}}
$$\ber=A_g\log{1\over(\e r)^n}+b_r,$$
and for each $q\in[1,+\infty]$ let us choose $r=r(\e,q)$  as follows:
$$
r^n={A_g\over2 C_1}\big(\log{1\over \e ^n}\big)^{1\over q'} \eqdef r
$$\def\aqn{{\a\over qn}}where $C_1$ is the constant in  in \eqref{T}.
Note that  for $\e$ small enough (independently of $q$)
$$\ber\ge \half A_g\log{1\over \e^n}\quad {\hbox{ if }} ,\qquad  {C_1r^n\over\ber}\le \Big(\log{1\over \e^n}\Big)^{-1/q}, \qquad  1\le q\le \infty.\eqdef{rrr}$$
Putting together the above  estimates we have that for $1\le q<\infty$
$$\eqalign{\big(\|\wtilde \phier\|_\na^{q\na}+\|T\wtilde \phier\|_\na^{q\na}\big)^{\aqn}&\le\Big((\ber+C)^q+C_1^qr^{nq}\Big)^{\aqn}\le \ber^\asn\bigg(\Big(1+{C\over \log{1\over\e^n}}\Big)^q+{1\over\log{1\over\e^n}}\bigg)^{\aqn}\cr&\le\ber^\asn\Big(1+{C\over\log{1\over \e^n}}\Big)\cr}$$
with  $C$ independent of $q$ (for example use the inequality $((1+x)^q+x)^{1/q}\le 1+2x$). 
For for $q=+\infty$ the same estimate holds, either by passing to the limit as $q\to+\infty $ (for fixed $\e$), or by direct check:
$$\eqalign{\max\big\{\|\wtilde \phier\|_\na,\|T\wtilde \phier\|_\na\big\}&\le \big(\max\{\ber+C,C_1r^n\}\big)^\asn\cr&\le\ber^{\asn}
\bigg(\max\Big\{1+{C\over\log{1\over\e^n}},1\Big\}\bigg)^\asn=\ber^\asn\Big(1+{C\over\log{1\over\e^n}}\Big).\cr}
$$

The value of the potential for $|x|\le \e r/2$ is estimated from below, using  \eqref{r4} and \eqref{TP}, as 
$$\eqalign{|T&\wtilde\phier(x)|\ge\bigg|\int_{\e r\le |y|\le r}\!\! |K(-y)|^{{\a\over n-\a}+1}+
\int_{\e r\le |y|\le r}\!\!\!\big(K(x-y)-K(-y)\big) K(-y)|K(-y)|^{{\a\over n-\a}-1}dy\bigg|\cr&-C
\ge \int_{\e r\le |y|\le r}|K(-y)|^\nna dy-\int_{\e r\le |y|\le r}|K(x-y)-K(-y)|\, |K(-y)|^{\a\over n-\a}dy-C\cr&\ge A_g \log{1\over(\e r)^n}-C+b_r -C\int_{|y|\ge\e r}
|x||y|^{-n-1}\ge \ber-C\ge \ber\Big(1-{C\over\log{1\over\e^n}}\Big),
\cr}\eqdef{12ggg}$$

If we now define\def\psier{\psi_{\e,r}} for $1\le q<\infty$
$$\psier={\wtilde\phier\over\Big(\|\wtilde \phier\|_\na^{q\na}+\|T\wtilde \phier\|_\na^{q\na}\Big)^{\a/qn}}, \eqdef{psier}$$
then $\|\psier\|_\na^{q\na}+\|T\psier\|_\na^{q\na}=1$, and using \eqref{rrr} we have, for $|x|\le \e r/2$
$$\eqalign{|T\psi_{\e,r}(x)|^\nna&\ge\bigg({ \ber\big(1-C/\log{1\over\e^n}\big)\over \ber^\asn\big(1+C/\log{1\over\e^n}\big)}\bigg)^{\nna}
\ge \ber\bigg(1-{C\over\log{1\over \e^n}}\bigg)\cr&= A_g\log{1\over( \e r)^n}+b_r\bigg(1-{C\over\log{1\over \e^n}}\bigg)-C.
\cr}$$
The same estimate follows in the case $q=+\infty$ if we define
$$\psier={\wtilde\phier\over\max\big\{\|\wtilde \phier\|_\na,\|T\wtilde \phier\|_\na\big\}}.$$
In the setting of Theorem 5 we have $q=1$, $r=1$, $b_r=0$ and  for any $\theta>1$
$$\eqalign{\int_{E\cap B_{\e /2}}\exp\bigg[{\theta\over A_g}\,|T \psi_{\epsilon,1}(x)|^{n\over n-\a}\bigg]dx
&\ge|E\cap B_{\e/2}|\exp\bigg[\theta\log{1\over\e^n}-C\bigg] \ge C\e^{-(\theta-1)n}\to+\infty\cr}$$
which  proves \eqref{12x} when $q=1$, and hence the sharpness of \eqref{12}.

In the setting of Theorem 6, we have $1<q\le+\infty$, and $K\notin L^\nna(|x|\ge1)$, which implies $b_r\to+\infty$. Hence we can estimate, with $r$ chosen as in \eqref r,
$$\eqalign{\int_{E\cap B_{\e r/2}}\!\!\!\exp\bigg[{1\over A_g}\,|T \psi_{\epsilon,r}(x)|^{n\over n-\a}\bigg]dx
&\ge|E\cap B_{\e r/2}|\exp\bigg[\log{1\over(\e r)^n}+{b_r\over A_g}\bigg(1-{C\over\log{1\over \e^n}}\bigg)-C\bigg]\cr& \ge C e^{{b_r\over2A_g}}\to+\infty,\cr}$$
as $\e\to0$, which proves \eqref{12x} when $q>1$, and hence the sharpness of \eqref{102}.

Still in the setting of Theorem 6 assume now $K$ homogeneous, in which case
$$b_r=\int_{1\le|y|\le r}|K(y)|^\nna dy=A_g\log r^n$$
and for any $\theta<1$ we take $r$ as in \eqref{r} and $\e$ as
$$\e^n=e^{-{1\over1-\theta}},\eqdef{teta}$$ and using the previous estimate we obtain 
$$\eqalign{\int_{E\cap B_{\e r/2}}&\exp\bigg[{\theta\over A_g}\,|T \psi_{\epsilon,r}(x)|^{n\over n-\a}\bigg]dx
\ge|E\cap B_{\e r/2}|\exp\bigg[\theta\log{1\over(\e r)^n}+{\theta b_r\over A_g}\bigg(1-{C\over\log{1\over \e^n}}\bigg)-C\bigg] \cr&\ge|E\cap B_{\e r/2}|\exp\bigg[\log{1\over\e^{n \theta}}-C\bigg]  \ge C r^n\e^{n(1-\theta)}= C(1-\theta)^{-{1\over q'}},\cr}$$
so that for fixed $E$ with positive measure, \eqref{102} with $C_\theta=C(1-\theta)^{-{1\over q'}}$ is reversed along the family $\{\psi_{\e,r}\}$ with $\e,r$  chosen as in \eqref r, \eqref{teta}, and $C$ depending possibly on $E$. Using Lemma~9 it's clear that \eqref{102aa} is also reversed along the same family.

Note that, if $K$ is homogeneous then 
 $$ \wtilde \phi_{\e,r}(y)=r^{-\a}\wtilde \phi_{\e,1}\Big({y\over r}\Big) $$
 and we can also use a dilation argument to arrive at the same  conclusions (see proof of Corollary 4).

The sharpness statements regarding \eqref{12aa} and \eqref{102aa} are obtained using the same extremal families as above, and using Lemma 9.

 In the vector-valued case the proof is completely similar. 
First, write an expansion as in \eqref{g} where each $p_j$ is a vector-valued polynomial whose components correspond to the Taylor formula of each $K_j(x-y)$. 
Then define
 $\phier$ as in \eqref{phie} 
and $\wtilde\phier$ as in \eqref{12g}, where the projection $P_{n}^r$ acts component-wise. The rest of the argument is exactly as in the scalar case. \endpf
\noin{\bf Remark.} The extremal family $\{\psi_{\e,1}\}$ in the previous proof can be used to show that the sharp constant in the original Adams inequality on bounded domains (or on domains of finite measure) cannot be improved under the more restrictive  Ruf condition.
\smallskip
\bigskip
\centerline{\bf 7. Proof of Theorem 7} 
\bigskip
We have that $T:D_0(T)\to L^\nsa$, where $D_0(T):=\{f\in L_c^{\na}(\R^n):\, Tf \in L^\nsa(\Rn)\},$  a subspace of $L^\nsa(\Rn)$.  The closability of $T$ follows at once from the next lemma:

\proclaim Lemma 16. If $\{f_k\}\subseteq D_0(T)$ is such that $f_k\buildrel{ L^{\nsa}}\over \longrightarrow 0$ and $T f_k\buildrel{ L^{\nsa}}\over \longrightarrow h$, then $h=0$ a.e.

\par 
\pf Proof. If $h$ is not zero a.e., we can assume that $\irn|h|^\nsa=1$
and that 
$$\int_{|x|\le R} |h|^\nsa\ge {3\over4},\qquad \int_{|x|\ge S} |h|^\nsa\le \e$$
for some $S>R>0$ and $\e$ small. Now consider
$$\varphi(x)={\rm sgn}(h)|h|^{n-\a\over\a}\chi_{B_R}^{}(x)$$
which is clearly in $L^{\nna}$ with $\|\varphi\|_\nna\le\|h\|_\na=1$, and has compact support, but its potential is not necessarily in $L^\nna$. For this to happen it is sufficient to normalize  $\varphi$ so that its mean  is zero, but we need to do this in a different way than the one used in the  proof of the sharpness statement of Theorem 5, which was localized inside a ball. 
\def\wvphi{{\widetilde\varphi}}

We let 
$$\wvphi(x)=\varphi(x)- \varphi(x-2S e_1)$$
and using \eqref{rl3}  we see that for $|x|\ge 4S$
$$\eqalign{|T\wvphi(x)|&\le \int_{|y|\le R}|K(x-y)-K(x-y+2Se_1)|\,| \varphi(y)|dy\cr&\le C_{\a,n}(2S)|x|^{\a-n-1}\int_{|y|\le R}|\varphi(y)|dy\le C|x|^{\a-n-1} \cr}$$
where $C$ depends on $R,S,\a,n$. Hence $T\wvphi$ is in $L^\nna$ for large $x$, and clearly this is also the case for small $x$.

Now we can say that $\wvphi\in L^\nna(\Rn)$ and $T\wvphi\in L^\nna(\Rn)$, and write
$$\eqalign{\irn& \wvphi h=\int_{B_R} |h|^\nsa-\int_{B_R}\varphi(x)h(x+2Se_1)dx\cr& 
\ge {3\over4}-\|\varphi\|_\nna\bigg(\int_{B_R}|h(x+2Se_1)|^\nsa dx\bigg)^\asn\ge {3\over4}-\e^{\asn}>{1\over2}\cr
}\eqdef{D1}$$
for $\e$ chosen small enough.

On the other hand, since $\wvphi,T\wvphi\in L^\nna$ and $f_k\to0,\,T f_k\to h$ in $L^\nsa$ we have 
$$\irn \wvphi h=\irn \wvphi \lim_{k\to\infty} T f_k=\lim_{k\to\infty} \irn \wvphi  T f_k=\lim_{k\to\infty} \irn (T \wvphi)  f_k=\irn (T\wvphi)   \lim_{k\to\infty} f_k=0$$
which contradicts \eqref{D1}\endpf


At this point we are in a position to apply a standard construction in order to close the operator $T$ see for ex. [Yo, Ch. II, Sect. 6]. Define
$$D(T)=\big\{f\in L^\nsa(\Rn): \,\exists \{f_k\}\subseteq D_0(T),\,\exists h\in L^\nsa(\Rn)\; {\rm with}  \;f_k\buildrel{ L^{\nsa}}\over \longrightarrow f,\;T f_k\buildrel{ L^{\nsa}}\over \longrightarrow h\big\}\eqdef{DD}$$
and because of Lemma 16 the function $h$ appearing in \eqref{DD} is independent of the sequence $f_k$, and the potential $Tf$ is well defined for $f$ in  $D(T)$, by letting $Tf=h$. The operator thus defined is the smallest closed extension of $T$ as defined on $D_0(T)$, and the class $D(T)$ is the closure of $D_0(T)$ under the graph norm 
$$\|f\|_{D(T)}^{}:=\big(\|f\|_\na^\na+\|Tf\|_\na^\na\big)^\asn.$$

The fact that Theorem 5 and Theorem 6 hold for $f\in D(T)$ is  now more or less straightforward. Indeed, given $f\in D(T)$, with  norm $\|f\|_{\na,q}:=\big(\|f\|_\na^{q \na}+\|Tf\|_\na^{q\na}\big)^{\a\over qn }\le 1$  there is  $f_k\in L_c^\na(\Rn)$ with $Tf_k\in L^\nsa(\Rn)$, and $f_k\to f$ in $D(T)$, i.e. $\|f_k-f\|_{\na,q}\to0$. All there is to do now is to write down the Adams inequalities in normalized form, i.e. using the functions $\wtilde f=f/\|f\|_{\na,q}\ge f,\, \wtilde f_k=f_k/\|f_k\|_{\na,q}\to  \wtilde f$ in $D(T)$, and the desired result follows from Fatou's lemma, (after possibly passing to a subsequence).

Now let us prove \eqref{sob1}. Let us first note that for $f\in D(I_\a)$ we have $\Da (c_\a I_\a f)=~f$. Indeed, if for some sequence  $\{f_k\}$ in $L^\nsa_c$ we have $f_k\to f$ and $I_\a f_k\to I_\a f$ in $L^\nsa$, then for all $\phi\in \S$ we also have
$$\eqalign{\int\Da (c_\a I_\a f) \phi&=\int (c_\a I_\a f)\Da\phi=\lim_{k\to\infty}\int (c_\a I_\a f_k)\Da \phi\cr&
=\lim_{k\to\infty}\int f_k (c_\a I_\a\Da \phi)=\lim_{k\to\infty}\int f_k  \phi=\int f\phi\cr}$$
where the first identity is due to the definition of $\Da$ on $L^p$ (as a tempered distribution), the second identity is due to H\"older's inequality ($\Da \phi\in L^p$ for any $p\ge 1$), the third identity is  due to Fubini's theorem, the fourth identity is true because $c_\a I_\a \Da \phi=\phi$ (take the FT of the left-hand side, which is in  $L^q$ for any  $q>\nna$, hence in $\S'$), and the fifth identity is again by H\"older's inequality.

Thus, $I_\a$ is injective on $D(I_\a)$, and if we define  temporarily 
$$\U_\a=\{u\in L^\nsa(\Rn): u=c_\a I_\a f,\;f\in D(I_\a)\},$$
endowed with the norm 
$$\|u\|_{\U_\a}^{}:=\Big(\|u\|_\na^\na+\|(-\Delta)^{\a\over2}u\|_\na^\na\Big)^{\a/n}=\|f\|_{D(I_\a)}^{},\qquad u=I_\a f,\, f\in D(I_\a)$$ then $I_\a:D(I_\a)\to \U_\a$ is a continuous bijection with inverse $\Da$, and $\U_\a$ is closed under such norm.
The identities 
$$(2\pi|x|)^{\a}=(1+\what h_1(x))(1+4\pi^2|x|^2)^{\a\over2},\;\;(1+4\pi^2|x|^2)^{\a\over2}=(1+\what h_2(x))\big(1+(2\pi|x|)^\a\big),$$
valid for some integrable functions $h_1,h_2$ (see [S, pp. 133-134]), combined with routine arguments show that 
$$W^{\a,\nsa}(\Rn)=\{u\in L^\nsa(\Rn):\, (-\Delta)^{\a\over2}u\in L^\nsa(\Rn)\}$$
and that  the norm $(\|u\|_\na^\na+\|(-\Delta)^{\a\over2}u\|_\na^\na)^{\a/n}$ is equivalent to $\|(I-\Delta)^{\a\over2}u\|_\na$ (this can actually be stated  for all $W^{\a,p}(\Rn)$, with $\a>0$ and $p>1$, see also [Hy1, Lemma A.3]).
Hence we have that $\U_\a\subseteq W^{\a,\nsa}$, and a Banach subspace. To prove the converse inclusion, by the density of  $C_c^\infty(\Rn)$ in $W^{\a,\nsa}(\Rn),$ we only need to prove that if $u\in C_c^\infty(\Rn)$ then there exists $\{f_k\}\in L_c^\nsa(\Rn)$ with 
$f_k\to f$ and $I_\a f_k\to u$, in $L^\nsa$, for some $f\in L^\nsa$. 

If $u\in C_c^\infty(\Rn)$ then there is $C>0$ (depending on $u$)  such that for all $h$ sufficiently large
$$|\Da u(x)|\le C|x|^{-\a-n},\qquad |x|\ge h.$$
This estimate is of course a special case of \eqref{schwarz}, but it can easily be proved as follows. Take $k\in \N$ and $\a_0\in(0,2)$ so that $\a=2k-\a_0$. Then we can write $\Da u=I_{\a_0}(-\Delta)^k u$, and the result follows easily after integrating by parts.

Also note that  $\Da u$ is orthogonal to any polynomial of degree less than $\a$, since the Fourier transform of such polynomial is a  linear combination of derivatives of the Dirac delta at $0$, having order strictly less than $\a$.

Our result will then be a consequence of the following lemma:

\proclaim Lemma 17. Let $f\in L^\nsa(\Rn)$ such that for some constant $C>0$ \smallskip
\item {i)} $\;|f(x)|\le C|x|^{-n-\a},\; |x|\ge h$\smallskip
\item{ii)} $f$ is orthogonal to  all polynomials with degree less than $\a$. \smallskip
Then, 
\smallskip \item{ a)} $I_\a f$ is well-defined a.e. and belongs to $L^\nsa(\Rn)$;\smallskip
\item{b)} If $f_k\in L^\nsa(\Rn)$ satisfies i) and ii) (with $C$ independent of $k$) and $f_k\to f$ in $L^\nsa$ then there is $\{f_{n_k}\}$ such that $I_\a f_{n_k}\to I_\a f$  in $L^\nsa$ ;
\smallskip\item{c)} There exist $f_k\in L^\nsa(\Rn)$ such that i) and ii) hold, $\supp f_k$ compact, and $f_k\to f$ in~$L^\nsa$.\par
\pf Proof of Lemma 17. a) Let us prove first that $I_\a f$ is finite a.e. and in $ L^\nsa(B_{2h})$, where $B_h=\{|x|\le h\}$. We write
$$I_\a f(x)= \int_{|y|\le h}|x-y|^{\a-n}f(y)dy+ \int_{|y|\ge h}|x-y|^{\a-n}f(y)dy:=J_1f(x)+J_2f(x).$$
Both integrals are clearly finite a.e. and Lemma 13 implies
$$\|J_1f\|_{L^\na(B_{2h})}\le C\|f\|_\na\eqdef{normaI1}$$
For $|x|\le 2h$
$$|J_2f(x)|\le J(x):= C\int_{|y|\ge h}|x-y|^{\a-n}|y|^{-\a-n}dy\le C\int_{|y|\ge h} |y|^{-2n}=C.$$
 Also for note that for $|x|\ge 2h$
$$J(x)\le C\int_{h\le|y|\le\half|x|}|x|^{\a-n}|y|^{-\a-n}dy+\int_{|y|\ge\half|x|}|y|^{-2n}dy\le C|x|^{-n}$$
so that $J\in L^\nsa(\Rn)$,
and in particular $I_\a f\in L^\nsa(B_{2h})$.
\smallskip Consider now the case $|x|\ge 2h$. From the Taylor's formula in \eqref{g}, we have, for all $m=1,2,....n$,
$$|x-y|^{\a-n}=\sum_{j=0}^{m-1}p_j(x,y)+R_m(x,y),\qquad |y|\le \half |x|\eqdef{g!}$$
where for $x\neq0$  $p_k(x,y)$ is a homogeneous polynomial of degree $k$ in the $y$ variable for $k\le m-1$, and such that estimates \eqref{stima} hold.
In particular,  
$$|R_m(x,y)|\le C|x|^{\a-n-m},\qquad |y|\le h,\, |x|\ge 2h.\eqdef{stimapm}$$

Now, let us write 
$$(0,n)=\Big(0,{n\over2}\Big)\cup\bigcup_{m=1}^n\Big[{n+m-1\over2},{n+m\over2}\Big)$$
so that for our given $\a$ either there is $m=1,....,n$ such that ${n+m-1\over2}\le\a<{n+m\over2}$, or else  $\a<{n\over2}$ in which case we let $m=0$. 
Since $f$ is orthogonal to polynomials of degree less than $\a$,  we can write, with $m$ chosen as above
$$\eqalign{I_\a f(x)&=\irn R_m(x,y)f(y)dy = \int_{|y|\le h}+\int_{|y|\ge h}:=\wtilde J_1f(x)+\wtilde J_2f(x).\cr} $$
Using \eqref{stimapm} (recall that we are assuming $|x|\ge 2h$)
$$|\wtilde J_1f(x)|\le C\int_{|y|\le h} |x|^{\a-n-m} |f(y)|dy\le C\|f\|_\na |x|^{\a-n-m}\in L^{\nsa}(|x|\ge 2h)\eqdef{stimaJ1}$$
and using \eqref{stima}
$$\eqalign{|\wtilde J_2f(x)|&\le \wtilde J(x):=C\int_{h\le |y|\le |x|/2}|R_m(x,y)||y|^{-\a-n}dy\cr&+C\int_{ |y|\ge |x|/2}|R_m(x,y)||y|^{-\a-n}
\le C\int_{h\le |y|\le |x|/2}|x|^{\a-n-m}|y|^{m-n-\a} dy+\cr&\hskip4em + C\int_{|y|\ge|x|/2}\bigg(|x-y|^{\a-n}+|x|^{\a-n}\sum_{j=0}^{m-1}|x|^{-j}|y|^j\bigg)|y|^{-\a-n}dy\cr& 
\hskip-3em =C|x|^{-n}+C|x|^{-n}\int_{|y|\ge\half}\bigg(|x^*-y|^{\a-n}+\sum_{j=0}^{m-1}|y|^j\bigg)|y|^{-\a-n}dy=C|x|^{-n}\in L^{\nsa}(|x|\ge 2h)\cr}\eqdef{stimaJ3}$$
(the function in $y$ in the last integral is integrable around $x^*$ and at infinity). This settles  a). 

To prove b) it is enough to show  that  if $f_k\to 0$ in $L^\nsa$, with $f_k\in L^\nsa(\Rn)$ satisfying i) and ii), then  up to a subsequence $I_\a f_k\to 0$ in $L^\nsa$. 
 With the above notation $I_\a f_k=J_1 f_k+J_2 f_k=\wtilde J_1f_k+\wtilde J_2 f_k$, a.e. We  will prove that up to a subsequence $J_1 f_k\to0$ and  $J_2 f_k\to0$ in $L^\nsa(B_{2h})$, as well as  $\wtilde J_1 f_k\to0$ and $\wtilde J_2 f_k\to0$ in $L^\nsa(|x|\ge 2h)$.
 
Using \eqref{normaI1}, \eqref{stimaJ1} we have
$$\|J_1f _k\|_{L^\na(B_{2h})}\le C\|f_k\|_\na\to0,\qquad\|\wtilde J_1 f_k\|_{L^\na(|x|\ge2h)}\le C\|f_k\|_\na\to 0$$
which show that $J_1f_k\to0$ in $L^\nsa(B_{2h})$ and $\wtilde J_1 f_k\to0$ in $L^\nsa(|x|\ge 2h)$.

On the other hand, 
$|J_2 f_k(x)|\le J(x)$ and  $|\wtilde J_2f_k(x)|\le \wtilde J(x)$ for a.e. $x$, and  with $J,\wtilde J\in L^\nsa(\Rn)$. Passing to a subsequence $f_{n_k}$ with $f_{n_k}\to 0$ a.e., we apply the dominated convergence theorem   to deduce first that $J_2f_{n_k}\to 0$, and  $\wtilde J_2f_{n_k}\to 0$ a.e., and then again to deduce that $J_2f_{n_k}\to0$ in $L^\nsa(B_{2h})$ and $\wtilde J_2 f_{n_k}\to0$ in $L^\nsa(|x|\ge 2h)$.  

To prove c), as in the proof of the sharpness statement of Theorem 5 consider  the space $\P_m$ of polynomials of degree up to $m$, with $m<\a\le m+1$, in the unit ball $B_1$, and with orthonormal basis $\{v_k\}_1^N$. Then $v_k^R(x)=R^{-n/2}v_k(x/R)$
defines an orthonormal basis on the ball $B_R=B(0,R)$, and we let
 $$f_R=f\chi_{B_R}^{},\qquad \wtilde f_R=f_R-\chi_{B_R}^{}\sum_{k=1}^N v_k^R\mu_k^R,\qquad \mu_k^R=\int_{B_R} v_k^R f.$$
Clearly $\wtilde  f_R$ is orthogonal to all polynomials of degree up to $m$ and we also have, using ii),
$$\eqalign{|\mu_k^R|&=\bigg|-\int_{|x|\ge R}R^{-n/2}v_k(x/R)f(x)dx\bigg|\le C R^{-n/2}\int_{|x|\ge R} |v_k(x/R)||x|^{-\a-n}dx\cr& =R^{-\a-n/2}\int_{|x|\ge1}|v_k(x)|x|^{-\a-n}dx\le CR^{-n/2-\a}\int_{|x|\ge1}|x|^{m-\a-n}=CR^{-n/2-\a}.\cr}$$
Hence, when $h\le |x|\le R$
$$|\wtilde f_R(x)|\le |f(x)|+C R^{-n-\a}\le C|x|^{-n-\a}$$
(since $|v_k^R(x)|\le CR^{-n/2}$ when $|x|\le R$), whereas $\wtilde f_R(x)=0$ when $|x|\ge R$. Finally, $f_R\to f$ in $L^\nsa$, and using the estimate above we have $|\wtilde f_R-f_R|\le C\chi_{B_R}^{}R^{-n-\a}\to0$ in $L^\nsa$. 

This concludes the proof of Lemma 17, and therefore of \eqref{sob1} (by taking $f=\Da u$),  and Theorem 7\endpf

\bigskip
\centerline{\bf 8. Proof of Theorem 3.}\bigskip  

In the case $P=\Da$ the inequalities \eqref{K4} and \eqref{K4a} are direct consequences of the corresponding inequalities in just Theorem 1 and Corollary 2, given the characterization \eqref{sob}, and the fact that the inverse of $\Da$ is precisely $c_\a I_\a$. 

In the remaining cases $\a$ is an integer, and we can assume that $u\in C_c^\infty(\Rn)$.

If $\a<n$ is an odd integer, writing $u=c_{\a+1}I_{\a+1}(-\Delta)^{\a+1\over2} u$ and integrating by parts gives
$$u(x)=J_\a f(x)=\irn c_{\a+1}(n-\a-1)|x-y|^{\a-n-1}(x-y)\cdot f(y)dy,\quad f=\nabla(-\Delta)^{\a-1\over2}u.\eqdef{J}$$
The kernel of  $J_\a$ changes sign component-wise, however we will verify that the alternate pointwise condition \eqref{12q} of   Theorem 5  holds, for our given $f$. Specifically, if $J_\a^+$ is the potential with kernel $c_{\a+1}(n-\a-1)|x-y|^{\a-n}$ and $f$ is as in \eqref{J}, 
then we can prove that for each $a\in \R^n$ 
$$|J_\a f(x)|=|u(x)|\le J_\a^+|f\chi_{|y-a|\le 2}^{} |(x)+C\|J_\a f\|_\na,\qquad|x-a|\le 1.\eqdef{Q}$$
It is enough to prove this for $a=0$ on the function $u_a(x)=u(x -a)$, so WLOG we can assume $a=0$.

Indeed, pick any smooth $\phi\in C_c^\infty(\Rn)$ such that $0\le|\phi|\le 1$ and 
$$\phi(y)=\cases {1 & if $|y|\le{3\over2}$\cr 0 & if $|y|\ge2$\cr}$$
and write, using Leibniz's rule and integration by parts (differentiations are in the $y$ variable),  
$$\eqalign{u(x)\phi(x)&=J_\a\big(\nabla\Delta^{\a-1\over2}(u\phi)\big)(x)=\int_{|y|\le2} c_{\a+1}\Delta|x-y|^{\a+1-n}\Delta^{\a-1\over2}(u\phi)(y)dy\cr&=
\int_{|y|\le2} c_{\a+1}\phi(y)\Delta|x-y|^{\a+1-n}\Delta^{\a-1\over2}u(y)dy+\cr&\hskip3em +
\int_{|y|\le2} c_{\a+1}\Delta|x-y|^{\a+1-n}\sum_{|k|+|h|\le {a-1}\atop |k|>0}b_{k,h,\a}(D^k\phi) (D^{h}u)\cr&
= -\int_{|y|\le2}  c_{\a+1}\phi(y)\nabla |x-y|^{\a+1-n}\cdot\nabla\Delta^{\a-1\over2}u(y)dy- \cr& \hskip3em 
-\int_{|y|\le2}  c_{\a+1}\Big(\nabla\phi(y)\cdot\nabla |x-y|^{\a+1-n}\Big)\Delta^{\a-1\over2}u(y)dy
+\cr&\hskip3em +\int_{|y|\le2} c_{\a+1}\sum_{|k|+|h|\le {a-1}\atop |k|>0}(-1)^{|h|}b_{k,h,\a}D_y^{h}\Big(\Delta|x-y|^{\a+1-n}D_y^k\phi(y)\Big)u(y)dy\cr}$$
 where $k=(k_1,...,k_n)$, $h=(h_1,...,h_n)$ are  multiindices, and the constants $b_{k,h,\a}$ are so that
$$\Delta^{{\a-1\over2}} (u\phi)=\sum_{|k|+|h|\le{\a-1}}b_{k,h,\a}(D^k\phi) (D^{h}u).$$
With further integrations by parts we can write
$$\eqalign{u(x)\phi(x)&
= -\int_{|y|\le2}  c_{\a+1}\phi(y)\nabla |x-y|^{\a+1-n}\cdot\nabla\Delta^{\a-1\over2}u(y)dy+ \cr& \hskip3em 
+ \sum_{0<|k|+|h|\le {\a+1}}c_{k,h,\a}\int_{|y|\le2} \Big(D_y^{h}|x-y|^{\a+1-n}D_y^k\phi(y)\Big)u(y)dy\cr}\eqdef{L} $$
for some other coefficients $c_{h,k,\a}$. Note that the derivatives of the function $\phi$ in the second term of \eqref{L}  all have positive order.
Now, for $|k|>0$ we have $\supp D^k\phi \subseteq \{{3\over2}\le |y|\le 2\}$, and for any fixed $x$ with $|x|\le 1$ the function  $y\to |x-y|^{\a+1-n}$ is $C^\infty$ outside the ball of radius ${3\over2}$, so that for all such $x$ 
$$\eqalign{|u(x)|&=|u(x)\phi(x)|\le\int_{|y|\le 2}c_{\a+1}\big|\nabla |x-y|^{\a+1-n}\big|\,\big|\nabla\Delta^{\a-1\over2}u(y)\big|dy+C\int_ {{3\over2}\le |y|\le2}|u(y)|dy\cr&\le \int_{|y|\le 2}c_{\a+1}(n-\a-1)| |x-y|^{\a-n}\big|\nabla\Delta^{\a-1\over2}u(y)\big|dy+C\|u\|_\na, \cr} $$
which is \eqref{Q}.

If $P$ is an elliptic operator as in the statement, the inequalities follow from Theorem~5, with $K=g_P$ as in \eqref{gp}, the homogeneous  fundamental solution of $P$. We cannot guarantee in general that $g_P$ does not change sign, however  it is possible to show that \eqref{12q} is true. Indeed, we have that  for $f=Pu$
$$|{g_P^{}}*f(x)|=|u(x)|\le |{g_P^{}}|*|f\chi_{|y-a|\le 2}^{}|(x)+C\|u\|_\na,\qquad  |x-a|\le1$$
the proof of which is a repetition of the proof of \eqref{Q}, but using the operator $P$ instead of $\nabla\Delta^{\a+1\over2}$.

To prove that the exponential constant in \eqref{K4}  is sharp when $P=\Da$,  it is enough to consider the functions $u_\e=c_\a I_\a\psi_\e $, where  the $\psi_\e$ were constructed in the proof of the sharpness statement in   Theorem 5 . In the case $\alpha$ odd, we can take the same extremal family $\{u_\e\}\in W_0^{\a,\nsa}(B(0,1))$ used in the original proof by Adams (see also  [FM1] proof of Theorem 6). Essentially,  $v_\e$   is a  smoothing of the function
$$\cases{0 & if $|y|\ge{3\over4}$\cr \log{1\over|y|^n} & if $2\e\le |y|\le {1\over2}$\cr
\log{1\over\e^n} & if $|y|\le \e$\cr}\eqdef{v}$$
which satisfies
$$\|v_\e\|_\na^\na\le C_2, \qquad
\|\nabla\Delta^{\a-1\over2}v_\e\|_\na^\na=
\gamma(P)^{n-\a\over\a}\log{1\over\e^n}+O(1),\eqdef{v1}$$ 
(some constant $C_2>0$) so that  the exponential integral in \eqref{K4} evaluated at the  functions 
$$u_\e:={v_\e\over (\|v_\e\|_\na^\na+\|\nabla\Delta^{\a-1\over2}v_\e\|_\na^\na)^{\asn}}\eqdef{v2}$$ can be made arbitrarily large if the exponential constant is larger than $\gamma(P)$.

To prove the sharpness of the exponential constant for $P$ elliptic, it is enough to take the family of functions $u_\e={g_P^{}}*\psi_\e$. 

The sharpness statements regarding \eqref{K4a} follow from the above results combined with Lemma 9.
\endpf
\noin{\bf Remark.}   In [IMM] the sharpness of the exponential constant $\gamma(P)=\pi$ in the case $P=(-\Delta)^{1\over4}$, $n=1$,  is proven in the form 
$$\sup_{\|u\|_2^2+\|(-\Delta)^{1\over4}u\|_2^2\le 1}\int_\sR |u|^a  \Big(e^{\pi u^2}-1\Big)dx=+\infty,\qquad a>2$$
and the authors ask whether the same statement holds for $0<a\le 2$. Our proof of Theorem~3  shows a lot more than this, namely, if $\Phi:[0,\infty)\to[0,\infty)$ is any measurable function such that $\lim_{t\to+\infty}\Phi(t)=+\infty$, then, in the setting of Theorem~3, and using the same extremal families,
$$\sup_{\|u\|_\na^\na+\|Pu\|_\na^\na\le 1}\int_{\sR^n}\Phi(|u|)\exp_{[\nsa-2]}\bigg[\gamma(P)|u(x)|^\nna\bigg]dx=+\infty\eqdef{sharper}$$
and a similar result hold for the general Adams inequalities in Theorem 5.
In the case $P=\Da$ the result above appears in a slightly stronger form in [Hy1, Thm.~1.1], where the extremal sequence in $W^{\a,\nsa}(\Rn)$  vanishes outside a given fixed ball. The result in [Hy1] is not stated on $\Rn$, but it can be easily deduced for that case, using the exponential regularization Lemma 9. 
\bigskip
\centerline{\bf 9. Proof of Corollary 4}
\bigskip
This proof is identical to the one of Theorem 6, in the case $K$ homogeneous. Given a function $u$ satisfying $\|u\|_\na^{qn/\a}+\|\nabla^\a u\|_\na^{q\na}\le 1$, we consider the functions $u_\lambda(x)=u(\lambda x)$, for $\lambda>0$, which, given the homogeneity of $P$, satisfy  
$$\|u_\lambda\|_\na^\na=\lambda^{-n}\|u\|_\na^\na,\qquad \|P u_\lambda\|_\na=\|P u\|_\na$$
and we choose $\lambda=\lambda(\theta)$ as we did in the proof of Theorem 6  to obtain inequality \eqref{102}.

The proof of the sharpness statement follows directly from the corresponding statement of Theorem 6 if $P$ is scalar, since we can just take $u_\theta=g_P*\psier$ with $\psier$ as in \eqref{psier} and $r,\e$ satisfying \eqref{teta}. If $P$ is not scalar, namely $P=\nabla(-\Delta)^{\a-1\over2}$, then we can use the following dilation argument, which works for the other $P$s as well. Take  the functions $v_\e,\,u_\e$ constructed in the previous proof,  which extremizes \eqref{K4}, and let $v_{\e,r}(x)=v_\e(x/r)$. Using \eqref{v1} and letting
$$r^n={\gamma(P)^{n-\a\over\a}\over C_2}\Big(\log{1\over\e^n}\Big)^{1\over q'},\qquad 1<q\le \infty$$
(with $C_2$ as in \eqref{v1}) we have
$$\eqalign{\|v_{\e,r}\|_{\na,q}&=\Big(r^{nq}\|v_\e\|_\na^{q\na}+\|Pv_\e\|_{\na}^{q\na}\Big)^{\a\over qn}\le
\Big(r^{nq}C_2^q+\big(\gamma(P)^{n-\a\over\a}\log{1\over\e^n}+C\big)^q\Big)^{\a\over qn}\cr&\le 
\gamma(P)^{n-\a\over n}\Big(\log{1\over \e^n}\Big)^{\asn}\bigg[\bigg({r^{n}C_2\over\gamma(P)^{n-\a\over\a}\log{1\over\e^n}}\bigg)^q+
\Big(1+{C\over\gamma(P)^{n-\a\over\a}\log{1\over\e^n}}\Big)^q\bigg]^{\a\over qn}\cr&
\le \gamma(P)^{n-\a\over n}\Big(\log{1\over \e^n}\Big)^\asn\Big(1+{C\over\log{1\over\e^n}}\Big),\cr}$$
and the same estimate holds for $\|v_{\e,r}\}_{\na,\infty}$. If we let 
$$u_{\e,r}= {v_{\e,r}\over\|v_{\e,r}\|_{\na,q}}  $$
then $\|u_{\e,r}\|_{\na,q}=1$, and when $|x|\le \e/2$
$$\bigg({|v_\e(x)|\over\|v_{\e,r}\|_{\na,q}}\bigg)^\nna={\big(\log{1\over\e^n}\big)^{\nna}\over \|v_{\e,r}\|_{\na,q}^\nna}\ge{1\over\gamma(P)}\log{1\over \e^n}-C$$
and the result follows from 
$$\eqalign{&\int_{E\cap B_{\e r/2}}
 \exp\bigg[{\theta \gamma(P)}|u_{\e,r}(x)|^\nna\bigg]dx=r^n \int_{r^{-1}(E\cap B_{\e r/2})}  
 \exp\bigg[{\theta\gamma(P)}\bigg({|v_\e(x)|\over\|v_{\e,r}\|_{\na,q}}\bigg)^\nna\bigg]dx\cr&
 \ge Cr^n\e^{n(1-\theta)}= C(1-\theta)^{-1\over q'}.\cr}$$
if we choose $\e^n=e^{-{1\over 1-\theta}}$.

\endpf

\bigskip
\centerline{\bf 10.  Inequalities for more general Borel measures}

\bigskip
The methods presented thus far allows us to obtain versions of the sharp inequalities in this paper when the non-regularized exponential is integrated against a positive Borel measure $\nu$ such that     
$$\nu\big(B(x,r)\big)\le Qr^{\sigma n},\qquad \forall x\in \Rn,\,\forall r>0\eqdef{bo}$$
for some $\sigma\in (0,1], \,Q>0$. However, to pass  from inequalities on sets of finite $\nu$ measure to the whole of $\Rn$, we cannot use  the exponential regularization Lemma 9  as is, since we are using two different measures in it. As  will be apparent from the proof below, we need to introduce some conditions at infinity satisfied by the measure $\nu$, in order to regularize the inequality on the whole space. It turns out that it is enough to ask that there are $r_1, Q'>0$ such that 
$$\nu(E)\le Q'|E|,\qquad \forall E {\hbox{ Borel  measurable with } } E\subseteq \{x: |x|\ge r_1\}.\eqdef{cond}$$
This condition is equivalent to asking that, outside a fixed ball, $\nu$ is  absolutely continuous with respect  to the Lebesgue measure, with bounded Radon-Nikodym derivative.
 An example is the singular measure $d\nu(x)=|x|^{(\sigma-1)n}dx$ considered in  [LL1], [LL2], [AY] and other papers.

\proclaim Theorem 18. Let $\nu$ be a positive Borel measure on $\Rn$ satisfying \eqref{bo}.  If  $K$ is a nonnegative or nonpositive Riesz-like kernel, then there exists a constant $C=C(\alpha,n,K,\sigma,Q)$ such that:
 \smallskip\item{(a)} For every measurable and compactly supported  $f:\Rn\to\R$ such that 
$$\|f\|_\na^\na+\|Tf\|_\na^\na\le 1\eqdef{bo1}$$
and for all Borel measurable $E\subseteq \R^n$ with $\nu(E)<\infty$, we have 
$$\int_E \exp\bigg[{\sigma\over A_g}|T f(x)|^\nna\bigg]d\nu(x)\le C(1+\nu(E)),\eqdef{bo2}$$
where $A_g$ is as in \eqref A. If in addition $\nu$  satisfies \eqref{cond} then 
$$\int_{\sR^n}\exp_{[\nsa-2]}\bigg[{\sigma\over A_g}|Tf(x)|^\nna\bigg]d\nu(x)\le C.\eqdef{bo3}$$
\smallskip\item{(b)} If $P$ is one of the operators  as in Theorem 3, then there is $C=C(\a,n,P,\sigma,Q)$ so that 
 for every  $u\in W^{\a,\nsa}(\Rn)$ with 
$$\|u\|_\na^\na+\|P u\|_\na^\na\le1\eqdef{bo4}$$
and for all Borel measurable $E\subseteq \Rn$ with $\nu(E)<\infty$, we have
$$\int_E\exp\Big[\sigma\gamma(P)|u(x)|^{n\over n-\a}\Big]d\nu(x)\le C(1+\nu(E)),\eqdef {bo5}$$
where $\gamma(P)$ is as in \eqref{gamma1}, \eqref{gamma2}. If in addition $\nu$ satisfies \eqref{cond} then 
$$\int_{\sR^n}\exp_{[\nsa-2]}\Big[\s\gamma(P)|u(x)|^{n\over n-\a}\Big]d\nu(x)\le C.\eqdef{bo6}$$
\item{(c)} If there exist $x_0,r_0$ such that 
$\nu( B(x_0,r))\ge c_1 r^{\sigma n}$, for $0<r<r_0$ with $c_1>0$ then the exponential constants in \eqref{bo3}, \eqref{bo6} are sharp. If  there exist $x_0,r_0$ such that 
$\nu(E\cap B(x_0,r))\ge c_1 r^{\sigma n}$, for $0<r<r_0$ with $c_1>0$, then the exponential constants  in the above inequalities are sharp.
\par

A word of caution:  the measure $\nu$ in this theorem enters only in the integration of the exponentials. The functions $f,\,I_\a*f,\, u,\,|\nabla^\a u|$ are still in $L^\nsa(\Rn)$ with respect to the Lebesgue measure.\medskip
\pf Proof. The proofs of inequalities \eqref{bo2}, \eqref{bo5} are identical to the corresponding  ones in Theorem 1 for the Lebesgue measure. The point is that the Adams inequality \eqref{12}  in   Theorem 5  holds for the measure $\nu$ as above, with the constant $A_g$ replaced by $\sigma^{-1}A_g$. The procedure is exactly the same, except now instead of \eqref{12h} we we use the following version of the  ``local  Adams inequality" 
$$\int_E \exp\bigg[{\sigma\over A_g}|T f(x)|^{n\over n-\alpha}\bigg]d\nu(x)\le C(1+|F|)\big(1+\log^+|F|+\nu(E)\big),\eqdef{localadams}$$
which can be obtained by tracking the constants carefully in the proof given in [FM1, Thm. 1] (See [FM4, Theorem 6]).
  Using \eqref{localadams} the entire proof given in  Theorem 5  goes through.

To deal with the \eqref{bo3}, \eqref{bo6} we need to modify Lemma 9. For simplicity we only prove \eqref{bo6} as the other inequality is completely similar. We assume condition \eqref{cond}
and estimate
$$\eqalign{\int_{\sR^n}&\exp_{[\nsa-2]}\Big[\s\gamma(P)|u(x)|^{n\over n-\a}\Big]d\nu(x)\le\cr& \le\int_{\{u\ge1\}}\exp\Big[\s\gamma(P)|u(x)|^{n\over n-\a}\Big]d\nu(x)+e^\alpha \int_{\{|u|\le 1\}}|u(x)|^{\nsa} d\nu(x)\cr}$$
  
now we have 
$$\eqalign{\nu\{|u|\ge1\}&= \nu\big(\{|u|\ge 1\}\cap B(0,r_1)\big)+ \nu\big(\{|u|\ge 1\}\cap B(0,r_1)^c\big)\cr&\le Qr_1^{\sigma n}+Q'|\{u\ge1\}|\le C(1+\|u\|_\na^\na)\cr}$$
so that we can use \eqref{bo5} to estimate the exponential integral over the set $\{|u|\ge 1\}$.
Finally, 
$$\eqalign{\int_{\{|u|\le 1\}}|u(x)|^{\nsa}d\nu(x)&\le \nu\big(B(0,r_1)\big)+\int_{\{|u|\le 1\}\cap\{|x|\ge r_1\}}|u(x)|^\nsa d\nu(x)\cr&\le Qr_1^{\sigma n}+Q'\int_{\{|u|\le 1\}}|u(x)|^\nsa dx\le C(1+\|u\|_\na^\na)\le C.\cr}$$
The proof of the sharpness result is the same as the one  in Theorem 5 and Theorem 3,   with the exception that this time $\nu(E\cap B_{\e/3})\ge c_1\e^{\sigma n}$ (using the notation in that proof, where we are taking $x_0=0$ and $\e<r_0$).

\endpf

\bigskip
\centerline{\bf 11. A sharp Trudinger  inequality on bounded domains}\centerline{\bf  without boundary conditions}
\bigskip
The next result has to do with smooth and  bounded domains, so it is somewhat unrelated to what we have done so far. We present it here since it is a nice and simple application of Lemma 10.

The original embedding result due to Trudinger [Tr, Thm 2] (see also [Po], [Yu]), states in particular that if $\Omega$ is open, bounded, and satisfies the cone condition, then for some constants $C,\gamma>0$
$$\int_\Omega e^{\gamma|u(x)|^{n\over n-1}}dx\le C,\qquad u\in W^{1,n}(\Omega),\;\| u\|_{W^{1,n}}\le 1\eqdef{tru}$$
where $\| u\|_{W^{1,n}}=\big(\|u\|_n^n+\|\nabla u\|_n^n\big)^{1/n}$ denotes the usual full norm in $W^{1,n}(\Omega)$. A version of this result for the space $W^{\a,\nsa}(\Omega)$ was given later by Strichartz [Str].

Sharp versions of \eqref{tru}    on smooth, connected, bounded domains $\Omega$ for functions $u\in~W^{\a,\nsa}(\Omega)$ are only known for $\alpha=1$ [CY], [Ci2], and for $\alpha=2$, if $\Omega$ is a ball [FM2].  In the case $\alpha=1$ Chang-Yang  ($n=2$) and Cianchi (any $n\ge2$) proved that there is $C$ such that
$$\int_\Omega \exp\bigg[2^{-{1\over n-1}}\gamma(\nabla)|u(x)-\overline u|^{n\over n-1}\bigg]dx\le C\qquad u\in W^{1,n}(\Omega),\;\|\nabla u\|_n\le 1,\eqdef{CC1}$$
where $\overline u=|\Omega|^{-1}\int_\Omega u$, and  $\gamma(\nabla)=n\omega_{n-1}^{1\over n-1}$ is the sharp constant for the Moser-Trudinger inequality on $W_0^{1,n}(\Omega)$, obtained in [Mo].

It is clear that some sort of normalization of $u$ is needed, as in \eqref{CC1}, if restrictions are imposed only on  the seminorm  $\|\nabla u\|_n$. Hence, it makes sense to ask about a sharp inequality under the full Sobolev norm condition $\|u\|_n^n+\|\nabla u\|_n^n\le 1$, and with no additional conditions on $u$, in the same spirit as in  the original paper by Trudinger. As far as we know no such result exists, however  we prove here that it can be  easily obtained by combining the Chang-Yang-Cianchi results and Lemma 10.

\proclaim Theorem 19. If $\Omega$ is a smooth, connected and bounded open set in $\Rn$, there exists a constant $C=C(\Omega)$ such that
$$\int_\Omega \exp\bigg[2^{-{1\over n-1}}\gamma(\nabla)|u(x)|^{n\over n-1}\bigg]dx\le C\eqdef{CC2}$$
 for each $u\in W^{1,n}(\Omega)$ with $\|u\|_n^n+\|\nabla u\|_n^n\le 1$. Moreover the exponential constant is sharp.\par
\pf Proof. In Lemma 10 take $\beta'=n$,  $V=Z=\big\{u\in W^{1,n}(\Omega): \ds\int_\Omega  u=0\big\}$, $\;T$ the identity on $V$, and $p(u)=\|\nabla u\|_n$. Then, using that 
$$|\overline u|\le |\Omega|^{-{1\over n}}\|u\|_n=|\Omega|^{-{1\over n}}\big(1-\|\nabla u\|_n^n\big)^{1\over n}$$
we obtain, using Chang-Yang-Cianchi's result together with (B) in Lemma 10
$$\eqalign{&\int_\Omega \exp\bigg[2^{-{1\over n-1}}\gamma(\nabla)|u(x)|^{n\over n-1}\bigg]dx\le\cr&\le\int_\Omega \exp\bigg[2^{-{1\over n-1}}\gamma(\nabla)\Big(|u(x)-\overline u|+|\Omega|^{-{1\over n}}\big(1-\|\nabla u\|_n^n\big)^{1\over n}\Big)^{n\over n-1}\bigg]dx\le Ce^{\gamma(\nabla)(2|\Omega|)^{-{1\over n-1}}}.\cr}$$
It is not hard to check that  if $0\in\p \Omega$ then the usual Moser sequence 
$$u_\e(x)=\cases{\log{1\over\e} & if $|x|<\e$\cr\log{1\over|x|} & if $\e\le|x|<1$\cr 0 & if  $|x|\ge1$\cr}$$
saturates the exponential constant in  \eqref{CC2}, arguing for example as in [F] pp. 451-453.

  The point is that as $\e\to0$ we have
 $$\|u_\e\|_{L^n(\Omega)}^n+\|\nabla u_\e\|_{L^n(\Omega)}^n\sim\half \|u_\e\|_{L^n(\sR^n)}^n+\half\| \nabla u_\e\|_{L^n(\sR^n)}^n\sim\half\| \nabla u_\e\|_{L^n(\sR^n)}^n= \half\omega_{n-1} \log{1\over\e}.$$
\endpf

We note that the results in [CY] and [Ci] were also obtained for smooth domains with finitely many conical singularities, in which case the sharp constant is $n(\theta_\Omega)^{1\over n-1}$, where $\theta_\Omega$ is the minimum solid aperture of the cones at the singularities. Needless to say, a result like Theorem 19 also holds under this more general situation, with the same sharp constant $n(\theta_\Omega)^{1\over n-1}$.
\eject
\bigskip\centerline{\bf Note added in proof}
\bigskip

Since this work was submitted,  a new, very interesting paper by Masmoudi-Sani has appeared [MS3]. In that paper the authors prove inequality (23) under condition (24) for the operators $\nabla^\a$, any $\alpha$ integer in $(0,n)$. Their proof is a clever combination of two tools: 1) a sharp growth estimate similar to the ones used by the authors in [27],  [29] but in the context  of second order borderline Lorentz-Sobolev spaces;   2)  a sharp embedding estimate for the borderline Sobolev space with Navier boundary conditions  into the second order borderline Lorentz-Sobolev space, which refines earlier work by Tarsi [Tar].\bigskip
\centerline{\bf References }
\vskip1em
\item{[AT]} Adachi S., Tanaka K., 
{\sl Trudinger type inequalities in $\R^N$ and their best exponents}, 
 Proc. Amer. Math. Soc. {\bf 128} (2000), 2051-2057.\smallskip 
\item{[A]} Adams D.R. {\sl
A sharp inequality of J. Moser for higher order derivatives},
Ann. of Math. {\bf128} (1988), no. 2, 385--398. 
\smallskip
\item{[AY]} Adimurthi, Yang Y., {\sl  An interpolation of Hardy inequality and Trundinger-Moser inequality in $\R^N$ and its applications}, IMRN {\bf 13} (2010), 2394-2426.\smallskip

\item{[Cao]} Cao D. M., 
{\sl Nontrivial solution of semilinear elliptic equation with critical exponent in $\R^2$}, 
Comm. Partial Differential Equations {\bf  17} (1992),  407-435. \smallskip
\item{[CST]} Cassani D., Sani F., Tarsi C.,
{\sl Equivalent Moser type inequalities in $\R^2$ and the zero mass case}, J. Funct. Anal. {\bf 267} (2014), 4236-4263.\smallskip 
\item{[CY]} Chang S.-Y.A., Yang  P.C., {\sl  Conformal deformation of metrics on $S^2$}, J. Differential Geom. {\bf 27} (1988), 259-296.\smallskip
\item{[Ci1]} Cianchi A., {\sl Moser-Trudinger trace inequalities}, Adv. Math. {\bf 217} (2008), 2005-2044.\smallskip
\item{[Ci2]} Cianchi A., {\sl Moser-Trudinger inequalities without boundary conditions and isoperimetric problems}, Indiana Univ. Math. J. {\bf 54} (2005), 669-705. \smallskip

\item{[CL]} Cohn W.S., Lu G., {\sl Best constants for Moser-Trudinger inequalities on the Heisenberg group}, Indiana Univ. Math. J. {\bf50} (2001), 1567-1591. \smallskip
\item{[do\'O]}  do \'O J.M.B. {\sl $N$-Laplacian equations in $\R^N$ with critical growth}, Abstr. Appl. Anal. {\bf 2} (1997), 301-315.\smallskip
\item{[F]}  Fontana L.,  {\sl Sharp borderline Sobolev inequalities on compact Riemannian manifolds}, Comment. Math. Helv. 
{\bf68} (1993), 415--454.\smallskip
\item{[FM1]} Fontana L., Morpurgo C., {\sl Adams inequalities on measure spaces}, Adv. Math. {\bf 226} (2011), 5066-5119.
\smallskip
\item{[FM2]} Fontana L., Morpurgo C., {\sl Sharp Moser-Trudinger inequalities for the Laplacian without boundary conditions}, J. Funct. Anal. {\bf 262} (2012),  2231-2271. 
\smallskip
\item{[FM3]} Fontana L., Morpurgo C., {\sl Critical integrability of fundamental solutions  and Adams inequalities  on spaces of infinite measure}, in preparation. \smallskip
\item{[FM4]} Fontana L., Morpurgo C., {\sl Sharp Adams and Moser-Trudinger inequalities on $\Rn$ and other spaces of infinite measure}, preprint (2015),	arXiv:1504.04678. \smallskip
\item{[GO]}  Grafakos L., Oh S., {\sl  The Kato-Ponce inequality},  Comm. Partial Differential Equations {\bf 39} (2014), 1128-1157.
\item{[Hy1]} Hyder A., {Moser functions and fractional Moser-Trudinger type
inequalities}, preprint (2015), arXiv:1510.06662.\smallskip
\item{[Hy2]} Hyder A., {\sl Structure of conformal metrics on $\Rn$ with constant $Q-$curvature},  preprint (2015), 	arXiv:1504.07095.\smallskip
\item{[IMN]} Ibrahim S., Masmoudi N., Nakanishi K., {\sl Trudinger-Moser inequality on the whole plane with the exact growth condition}, J. Eur. Math. Soc. (JEMS) {\bf 17} (2015), 819-835.\smallskip
\item{[IMM]} Iula S., Maalaoui A., Martinazzi L., {A fractional Moser-Trudinger type inequality in one dimension and its critical points}, Differential Integral Equations {\bf  29} (2016), 455-492. \smallskip
\item{[L]} Lieb E.H., {\sl Sharp constants in the 
Hardy-Littlewood-Sobolev and related inequalities}, 
  Ann. of Math. {\bf118} (1983), 349-374.\smallskip
\item{[LL1]} Lam N., Lu G., {\sl 
Sharp Moser-Trudinger inequality on the Heisenberg group at the critical case and applications}, Adv. Math. {\bf 231} (2012),  3259-3287.\smallskip 
\item{[LL2]} Lam N., Lu G., {\sl A new approach to sharp Moser-Trudinger and Adams type inequalities: a rearrangement-free argument}, 
J. Differential Equations {\bf 255} (2013), 298-325.\smallskip 
\item{[LLZ]} Lam N., Lu G., Zhang L., {\sl Equivalence of critical and subcritical sharp Trudinger-Moser-Adams inequalities},  arXiv:1504.04858 (2015).\smallskip
\item{[LR]} Li Y., Ruf B., {\sl A sharp Trudinger-Moser type inequality for unbounded domains in $R^n$},   
 Indiana Univ. Math. J. {\bf 57} (2008),  451-480.\smallskip

\item{[LTZ]} Lu G., Tang  H., Zhu M., {\sl Best constants for Adams' inequalities with the exact growth condition in $\R^n$}, Adv. Nonlinear Stud. {\bf 15} (2015), 763-788.\smallskip

\item{[MS1]} Masmoudi N.; Sani F.,{ \sl Trudinger-Moser inequalities with the exact growth condition in $\R^n$ and applications}, Comm. Partial Differential Equations {\bf 40} (2015), 1408-1440.\smallskip
\item{[MS2]}  Masmoudi N., Sani F., {\sl Adams' inequality with the exact growth condition in $\R^4$}, Comm. Pure Appl. Math. {\bf 67} (2014), 1307-1335.\smallskip
\item{[MS3]}  Masmoudi N., Sani F., {\sl Higher order Adams' inequality with the exact growth condition}, Commun. Contemp. Math. (2017).\smallskip
\item{[Mo]} Moser J., {\sl A sharp form of an inequality by N. Trudinger}, Indiana Univ. Math. J. {\bf20} (1970/71), 1077-1092.\smallskip

\item{[Og]} Ogawa T. {\sl 
A proof of Trudinger's inequality and its application to nonlinear Schrödinger equations}, 
Nonlinear Anal. {\bf 14} (1990),  765-769.\smallskip 

\item{[Oz]}  Ozawa T., {\sl On critical cases of Sobolev's inequalities}, J. Funct. Anal. {\bf 127} (1995), 259-269.\smallskip 
\item{[Pa]} Panda R. {\sl 
Nontrivial solution of a quasilinear elliptic equation with critical growth in $\R^n$}, Proc. Indian Acad. Sci. Math. Sci. {\bf 105} (1995), 425-444.\smallskip 
\item{[Po]} Pohozhaev  S.I., {\sl On the imbedding Sobolev theorem for pl = n}, Doklady Conference,
Section Math. Moscow Power Inst. (1965), 158-170 (Russian).\smallskip
\item{[Ruf]} Ruf B., {\sl A sharp Trudinger--Moser type inequality for unbounded domains in $\R^2$}, J. Funct. Anal. {\bf 219} (2005), 340-367.\smallskip
\item{[RS]} Ruf B., Sani F., {\sl
 Sharp Adams-type inequalities in $\R^n$},  
Trans. Amer. Math. Soc. {\bf 365} (2013), 645-670.\smallskip 
\item{[Sa]} Samko S.G., {\sl Hypersingular integrals and their applications}, Analytical Methods and Special Functions, 5. Taylor \& Francis, Ltd., London, 2002.\smallskip
\item{[S]} Stein E.M., {\sl
Singular integrals and differentiability properties of functions}, 
Princeton Mathematical Series {\bf 30}, Princeton University Press, 1970.\smallskip
\item{[Str]} Strichartz R.S., {\sl A note on Trudinger's extension of Sobolev's inequalities}, Indiana Univ. Math. J. {\bf 21} (1971/72), 841-842.\smallskip
\item{[Tar]} Tarsi C., 
{\sl Adams' inequality and limiting Sobolev embeddings into Zygmund spaces}, 
Potential Anal. {\bf 37} (2012),  353-385.\smallskip 

\item{[Tr]} Trudinger N.S., {\sl 
On imbeddings into Orlicz spaces and some applications}, 
J. Math. Mech. {\bf17} (1967), 473-483. 
\smallskip
\item{[Y]} Yang Y., {\sl 
Trudinger-Moser inequalities on complete noncompact Riemannian manifolds},  
J. Funct. Anal. {\bf 263} (2012),  1894-1938.\smallskip 
\item{[Yo]} Yosida K., {\sl  Functional Analysis}, Sixth edition. Springer-Verlag, Berlin-New York, 1980.

\item{[Yu]} Yudovic V.I., {\sl Some estimates connected with integral operators and with solutions of
elliptic equations}, Dokl. Akad. Nauk SSSR {\bf138} (1961), 804-808, English translation in
Soviet Math. Doklady 2 (1961), 746-749.\smallskip
\smallskip\smallskip
\noin Luigi Fontana \hskip19em Carlo Morpurgo

\noin Dipartimento di Matematica ed Applicazioni \hskip5.5em Department of Mathematics 
 
\noin Universit\'a di Milano-Bicocca\hskip 12.8em University of Missouri

\noin Via Cozzi, 53 \hskip 19.3em Columbia, Missouri 65211

\noin 20125 Milano - Italy\hskip 16.6em USA 
\smallskip\noin luigi.fontana@unimib.it\hskip 15.3em morpurgoc@missouri.edu

\end